\documentclass{amsart}

\usepackage[mathcal]{eucal}
\usepackage{amsmath, amssymb, graphics, labelfig}

\newtheorem{thm}{Theorem}
\newtheorem{prop}[thm]
{Proposition}
\newtheorem{lem}[thm]{Lemma}

\theoremstyle{remark}

\newcommand{\col}{\kern -3pt :}
\renewcommand{\leq}{\leqslant}
\renewcommand{\geq}{\geqslant}

\begin{document}

\title[Quantum Teichm\"uller
space and surface
diffeomorphisms]
{Representations
of
the quantum Teichm\"uller
space and invariants of
surface diffeomorphisms}

%\authors{Francis Bonahon\\ Xiaobo Liu}
%\address{Department of
%Mathematics, University of
%Southern California, Los Angeles,
%CA~90089-2532, U.S.A.}

%\secondaddress{Department of
%Mathematics,  Columbia University,
%2990 Broadway, New York, NY~10027, 
%U.S.A.}

\author{Francis Bonahon}
\address{Department of
Mathematics, University of
Southern California, Los Angeles,
CA~90089-2532, U.S.A.}
\author{Xiaobo Liu}
\address{Department of
Mathematics,  Columbia University,
2990 Broadway, New York, NY~10027, 
U.S.A.}

\begin{abstract} We investigate
the representation theory of the
polynomial core $\mathcal T_S^q$
of the quantum Teichm\"uller space
of a punctured surface
$S$. This is a purely algebraic
object, closely related to the
combinatorics of the simplicial
complex of ideal cell
decompositions of
$S$. Our main result is that
irreducible finite-dimensional
representations of  $\mathcal
T_S^q$ are classified, up to
finitely many choices, by group
homomorphisms from the
fundamental group
$\pi_1(S)$ to the isometry group
of the hyperbolic 3--space
$\mathbb H^3$. We exploit this
connection between algebra and
hyperbolic geometry to exhibit
invariants of diffeomorphisms of
$S$. 
\end{abstract}

%\primaryclass{57R56}                
%\secondaryclass{57M50, 20G42} 

\maketitle
%\maketitlepage

This work finds its motivation in
the emergence of various
conjectural connections between
topological quantum field theory
and hyperbolic geometry, such as
 the now famous Volume
Conjecture of Rinat Kashaev
\cite{Kash2}, Hitochi Murakami and
Jun Murakami
\cite{MM}. For a hyperbolic
link $L$ in the 3--sphere
$S^3$, this conjectures relates the hyperbolic
volume of the complement $S^3-L$
to the asymptotic behavior of the
$N$--th colored Jones polynomial
$J^N_L(\mathrm e^{2\pi \mathrm
i/N})$ of $L$, evaluated at the
primitive
$N$--th root of unity $e^{2\pi
\mathrm i/N}$. At this point, the heuristic
evidence
\cite{Kash2, Mur, Yok2, Yok}  for
the Volume Conjecture is based on
the observation
\cite{Kash2, MM} that the
$N$--th Jones polynomial can be
computed using an explicit
R-matrix whose asymptotic
behavior is related to Euler's
dilogarithm function, which is
well-known to give the hyperbolic
volume of an ideal tetrahedron in
$\mathbb H^3$ in terms of the
cross-ratio of its vertices. We
wanted to establish a more
conceptual connection between the
two points of view, namely between
quantum algebra and 3--dimensional
hyperbolic geometry. 

We investigate such a
relationship, 
provided by the quantization of the
Teich\-m\"ul\-ler space of a surface,
as developed by Rinat Kashaev
\cite{Kash3}, Leonid Chekhov and
Vladimir Fock \cite{CF}. More
precisely, we follow the
exponential version of the
Chekhov-Fock approach. This
enables us to formulate our
discussion in terms of
non-commutative algebraic
geometry and finite-dimensional
representations of algebras,
instead of Lie algebras and
self-adjoint operators of Hilbert
spaces. This may be physically
less relevant, but this point of
view is better adapted to the
problems that  we have in mind. The
mathematical foundations of this
non-commutative algebraic
geometric point of view are
rigorously established in
\cite{Liu}.

More precisely, let $S$ be a
surface of finite topological
type, with genus $g$ and with
$p\geq 1$ punctures. An
\emph{ideal triangulation} of $S$
is a proper 1--dimensional
submanifold  whose complementary
regions are infinite triangles
with vertices at infinity, namely
at the punctures. For an ideal
triangulation $\lambda$ and a
number $q= \mathrm e ^{\pi
\mathrm i
\hbar} \in \mathbb C$, the
\emph{Chekhov-Fock algebra}
$\mathcal T_\lambda^q$ is the
algebra over $\mathbb C$ defined
by   generators
$X_1^{\pm1}$, $X_2^{\pm1}$,
\dots, $X_n^{\pm1}$ associated to
the components of $\lambda$ and
by relations
$X_iX_j=q^{2\sigma_{ij}}X_jX_i$,
where the $\sigma_{ij}$ are
integers determined by the
combinatorics of the ideal
triangulation $\lambda$. This
algebra has a well-defined
fraction division algebra
$\widehat{\mathcal T}_\lambda^q$.
In concrete terms, 
$\mathcal T_\lambda^q$ consists
of the formal Laurent polynomials
in  variables $X_i$ satisfying
the skew-commutativity relations
$X_iX_j=q^{2\sigma_{ij}}X_jX_i$,
while its fraction algebra
$\widehat{\mathcal T}_\lambda^q$
consists of formal rational
fractions in the $X_i$ satisfying
the same relations. 

As one moves from one ideal
triangulation
$\lambda$  to another 
$\lambda'$, Chekhov and Fock
\cite{Fo, FC, CF} (see also
\cite{Liu}) introduce
\emph{coordinate change
isomorphisms} 
$\Phi_{\lambda\lambda'}^q \col
\widehat{\mathcal T}_{\lambda'}
^q \rightarrow
\widehat{\mathcal T}_\lambda^q$
which satisfy the natural
property that
$\Phi_{\lambda''\lambda'}^q \circ
\Phi_{\lambda'\lambda}^q =
\Phi_{\lambda''\lambda}^q $ for
every ideal triangulations
$\lambda$, $\lambda'$,
$\lambda''$. In a triangulation
independent way, this associates
to the surface
$S$ the algebra
$\widehat{\mathcal T}_S^q$ defined
as  the quotient of the family of
all
$\widehat{\mathcal T}_\lambda^q$,
with $\lambda$ ranging over ideal
triangulations of the surface
$S$, by the equivalence relation
that identifies 
$\widehat{\mathcal T}_\lambda^q$
and
$\widehat{\mathcal T}_{\lambda'}
^q $ by the coordinate change
isomorphism
$\Phi_{\lambda\lambda'}^q$. By
definition, $\widehat{\mathcal
T}_S^q$ is the \emph{quantum
Teichm\"uller space} of the
surface $S$.

This construction and definition
are motivated by the case where
$q=1$, in which case
$\widehat{\mathcal T}_\lambda^1$
is just the algebra
$\mathbb C (X_1, X_2, \dots,
X_n)$ of rational functions in $n$
commuting variables. Bill
Thurston associated to each ideal
triangulation a global coordinate
system for the
\emph{Teichm\"uller space}
$\mathcal T(S)$ consisting of all
isotopy classes of complete
hyperbolic metrics on
$S$. Given two ideal
triangulations $\lambda$ and
$\lambda'$, the corresponding
coordinate changes are rational,
so that there is a well-defined
notion of rational functions on
$\mathcal T(S)$. For a given
ideal triangulation
$\lambda$, Thurston's shear
coordinates provide a canonical
isomorphism between the algebra
of rational functions on
$\mathcal T(S)$ and  $\mathbb C
(X_1, X_2, \dots, X_n) \cong 
\widehat{\mathcal T}_\lambda^1$.
It turns out that the
$\Phi_{\lambda\lambda'}^1$ are
just the corresponding coordinate
changes. Therefore, the quantum
Teichm\"uller space
$\widehat{\mathcal T}_S^q$ is a
(non-commutative) deformation of
the algebra of rational functions
on the Teichm\"uller space
$\mathcal T(S)$.

Although the construction of
$\widehat{\mathcal T}_S^q$ was
motivated by the geometry, a
result of Hua Bai \cite{Bai}
shows that it actually depends
only on the combinatorics of ideal
triangulations. Indeed, once we
fix the definition of the
Chekhov-Fock algebras ${\mathcal
T}_\lambda^q$, the coordinate
change isomorphisms 
$\Phi_{\lambda\lambda'}^q \col
\widehat{\mathcal T}_{\lambda'}
^q \rightarrow
\widehat{\mathcal T}_\lambda^q$
are uniquely determined if we
require them to satisfy a certain
number of natural conditions, a
typical one being the locality
condition: if
$\lambda$ and $\lambda'$ share a
component
$\lambda_i$ as well as any
component of $\lambda$ that is
adjacent to $\lambda_i$, then
$\Phi_{\lambda\lambda'}^q$ must
respect the corresponding
generator 
$X_i$. 

A standard method to move from abstract algebraic constructions to more concrete applications is to consider finite-dimensional representations. In the case of algebras, this means 
algebra homomorphisms valued in
the algebra $\mathrm{End}(V)$ of
endomorphisms of a
finite-dimensional vector space
$V$ over $\mathbb C$. Elementary
considerations show that these
can exist only when $q$ is a root
of unity.  

\begin{thm}
\label{thm:RepCheFockIntro}
 Suppose that
$q^2$ is a primitive $N$--th root
of unity, and consider the
Chekhov-Fock algebra $ \mathcal
T^q_\lambda$ associated to an
ideal triangulation
$\lambda$. Every irreducible
finite-dimensional  representation
of $ \mathcal T^q_\lambda$ has
dimension
$N^{3g+p-3}$ if $N$ is odd, and
$N^{3g+p-3}/2^g$ if $N$ is even,
where $g$ is the genus of the
surface $S$ and where $p$ is its
number of punctures. Up to
isomorphism, such a
representation is classified by:
\begin{enumerate}
\item a non-zero complex number
$x_i \in \mathbb C^*$ associated
to each edge of $\lambda$;
\item a choice of an $N$--th root
for each of
$p$ explicit monomials in the
numbers $x_i$;
\item when $N$ is even, a choice
of square root for each of $2g$
explicit monomials in the numbers
$x_i$. 
\end{enumerate} Conversely, any
such data can be realized by an
irreducible finite-dimensional
representation of $ \mathcal
T^q_\lambda$.
\end{thm}

The numbers $x_i \in
\mathbb C^*$ appearing in the
classification of a
representation $\rho\col
\mathcal T_\lambda^q \rightarrow
\mathrm{End}(V)$ are characterized
by the property that
$\rho(X_i^N) = x_i\,
\mathrm{Id}_V$ for the
corresponding generator
$X_i$ of $\mathcal
T_\lambda^q$. Theorem~\ref
{thm:RepCheFockIntro} is proved in
Section~\ref{sect:RepTheory}.   The
main step in the proof, which has a
strong topological component, is to
determine the algebraic structure
of the algebra $
\mathcal T^q_\lambda$ and is
completed in
Section~\ref{sect:StructureCF}
after preliminary work in 
Section~\ref{sect:StructureWP}.
Another important feature of
Theorem~\ref
{thm:RepCheFockIntro}  is the way
it is stated, which closely ties
the classification to the
combinatorics of the ideal
triangulation
$\lambda$ in
$S$ and counterbalances the fact
that the structure results for 
$\mathcal T^q_\lambda$ are not very
explicit. 

Theorem~\ref
{thm:RepCheFockIntro}  shows that
the Chekhov-Fock algebra has a
rich representation theory.
Unfortunately, for dimension
reasons, its fraction algebra
$\widehat{\mathcal T}^q_\lambda$
and, consequently, the quantum
Teichm\"uller space
$\widehat{\mathcal T}^q_S$ cannot
have any finite-dimensional
representation. This leads us to
introduce the
\emph{polynomial core} $\mathcal
T_S^q$ of the quantum
Teichm\"uller space
$\widehat{\mathcal T}^q_S$,
defined as the family $\{ \mathcal
T_\lambda^q\} _{\lambda \in
\Lambda(S)}$ of all Chekhov-Fock
algebras $\mathcal T_\lambda^q$,
considered as subalgebras of
$\widehat{\mathcal T}^q_S$, as
$\lambda$ ranges over the set
$\Lambda(S)$ of all isotopy
classes of ideal triangulations
of the surface
$S$.  In
Section~\ref{sect:PolCore}, we
introduce and analyze the
consistency of a notion of
representation of the polynomial
core, consisting of  the data of
 representations
$\rho_\lambda\col
\mathcal T^q_\lambda
\rightarrow
\mathrm{End} (V)$ for all 
$\lambda \in \Lambda(S)$ that
behave well under the coordinate
changes $
\Phi^q_{\lambda\lambda'} $.

We now jump from the purely
algebraic representation theory of
the polynomial core  $\mathcal
T_S^q$ to 3--dimensional
hyperbolic geometry.  Theorem~\ref
{thm:RepCheFockIntro} says that,
up to a finite number of choices,
an irreducible representation of 
$\mathcal T^q_\lambda$ is
classified by certain numbers
$x_i\in \mathbb C^*$ associated
to the edges of the ideal
triangulation 
$\lambda$ of $S$. There is a
classical geometric object which is
also associated to $\lambda$ with
the same edge weights $x_i$.
Namely, we can consider in the hyperbolic
3--space
$\mathbb H^3$ the
pleated surface that has pleating
locus
$\lambda$, that has shear
parameter along the
$i$--th edge of $\lambda$ equal to
the real part of $\log x_i$, and
that has bending angle along this
edge equal to the imaginary part
of $\log x_i$. In turn, this
pleated surface has a
\emph{monodromy representation},
namely a group homomorphism from
the fundamental group $\pi_1(S)$
to the  group
$\mathrm{Isom}^+(\mathbb H^3)
\cong \mathrm{PSL}_2(\mathbb C)$
of orientation-preserving
isometries of $\mathbb H^3$. This
construction associates to a
representation of the Chekhov-Fock
algebra $ \mathcal T^q_\lambda$ a
group homomorphism
$r\col \pi_1(S) \rightarrow
\mathrm{PSL}_2(\mathbb C)$,
well-defined up to conjugation by
an element of
$\mathrm{PSL}_2(\mathbb C)$.

It turns out that, for a suitable
choice of $q$, this construction
is well-behaved under coordinate
changes.  The fact that $q^2$ is
a primitive $N$--th root of unity
implies that
$q^N=\pm1$, but the following
result  requires that
$q^N=(-1)^{N+1}$. This is
automatically satisfied if $N$ is
even. 

\begin{thm} 
\label{thm:HyperShadowIntro} Let
$q$ be a primitive $N$--th root
of $(-1)^{N+1}$, for instance
$q=-\mathrm e^{2 \pi
\mathrm i/N}$. If 
$\rho =
\{\rho_\lambda\col
\mathcal T^q_\lambda
\rightarrow
\mathrm{End} (V)\}_{\lambda \in
\Lambda(S)}$ is a
finite-dimensional irreducible
representation of the polynomial
core $\mathcal T^q_S$ of the
quantum Teichm\"uller space
$\widehat{\mathcal T}^q_S$, the
representations $\rho_\lambda$
induce the same monodromy
homomorphism $r_\rho\col \pi_1(S)
\rightarrow
\mathrm{PSL}_2(\mathbb C)$.
\end{thm}

Theorem~\ref{thm:HyperShadowIntro}
is essentially equivalent to the
property that, for the choice of
$q$ indicated, the pleated
surfaces respectively associated
to the representations 
$\rho_\lambda\col
\mathcal T^q_\lambda
\rightarrow
\mathrm{End} (V)$ and 
$\rho_\lambda \circ
\Phi_{\lambda\lambda'}^q\col
\mathcal T^q_{\lambda'}
\rightarrow
\mathrm{End} (V)$ have (different
pleating loci but) the same
monodromy representation 
$r_\rho\col \pi_1(S)
\rightarrow
\mathrm{PSL}_2(\mathbb C)$. Its
proof splits into two parts: a
purely algebraic computation in
Section~\ref
{sect:NonQuantumShadow}, which is
based on the quantum binomial
formula and is borrowed from a
remark in
\cite{FC}, relates the quantum case
to the non-quantum case where
$q=1$;  a more geometric part in
Section~\ref{sect:HyperShadow}
 is completely centered on
the non-quantum situation and uses
pleated surfaces in hyperbolic
3--space.

The homomorphism $r_\rho$ is the
\emph{hyperbolic shadow} of the
representation $\rho$. Not every
homomorphism $r\col \pi_1(S)
\rightarrow
\mathrm{PSL}_2(\mathbb C)$ is the
hyperbolic shadow of a
representation of the polynomial
core, but many of them are:

\begin{thm}
\label{thm:HyperGivesRepIntro} An
injective homomorphism  
$r\col \pi_1(S)
\rightarrow
\mathrm{PSL}_2(\mathbb C)$ is the
hyperbolic shadow of a finite
number of irreducible
finite-dimensional
representations of the polynomial
core $\mathcal T^q_S$, up to
isomorphism. More precisely, this
number of representations is
equal to $2^lN^p$ if $N$ is odd,
and $ 2^{2g+l}N^p$ if $N$ is
even, where $g$ is the genus of
$S$, $p$ is its number of
punctures, and $l$ is the number
of ends of $S$ whose image under
$r$ is loxodromic.
\end{thm}

As an application of this
machinery, we construct new and
still mysterious invariants of
(isotopy classes of) surface
diffeomorphisms, by using
Theorems~\ref
{thm:HyperShadowIntro} and
\ref{thm:HyperGivesRepIntro} to
go back and forth between
hyperbolic geometry and
representations of the polynomial
core $\mathcal T^q_S$.

Let $\varphi$ be a diffeomorphism
of the surface $S$. Suppose in
addition that $\varphi$ is
homotopically aperiodic
(also called homotopically
pseudo-Anosov), so that its
(3--dimensional) mapping torus
$M_\varphi$ admits a complete
hyperbolic metric. The hyperbolic
metric of $M_\varphi$ gives an
injective homomorphism 
$r_\varphi\col \pi_1(S)
\rightarrow
\mathrm{PSL}_2(\mathbb C)$ such
that $r_\varphi \circ \varphi^*$
is conjugate to $r_\varphi$, where
$\varphi^*$ is the isomorphism of
$\pi_1(S)$ induced by $\varphi$.

The diffeomorphism $\varphi$ also
acts on the quantum Teichm\"uller
space and on its polynomial core $
\mathcal T_S^q$. In particular,
it acts on the set of
representations of
$\mathcal T_S^q$ and, because
$r_\varphi
\circ
\varphi^*$ is conjugate to
$r_\varphi$, it sends a
representation with hyperbolic
shadow $r_\varphi$ to another
representation with shadow
$r_\varphi$. 
Actually, when $N$ is odd,
there is a preferred
representation 
$\rho_\varphi$ of $\mathcal
T_S^q$ which
is fixed by the action of
$\varphi$, up to isomorphism. 
This statement means that, for
every ideal triangulation
$\lambda$, we have a
representation
$\rho_\lambda\col
\mathcal T^q_\lambda
\rightarrow
\mathrm{End} (V)$ of dimension 
$N^{3g+p-3}$ and an isomorphism
$L_\varphi^q$ of $V$ such that
\begin{equation*}
\rho_{\varphi(\lambda)} \circ
\Phi_{\varphi(\lambda) \lambda}
(X) = L_\varphi^q \cdot
\rho_\lambda (X)
\cdot (L_\varphi^q)^{-1}
\end{equation*} in
$\mathrm{End}(V)$ for every
$X\in \mathcal T_\lambda^q$, for
a suitable interpretation of the
left hand side of the equation.

\begin{thm} 
\label{thm:intro:Invariant}
Let $N$ be odd. Up to
conjugation and up to
multiplication by a constant, the
isomorphism
$L_\varphi^q$ depends only on the
homotopically aperiodic
diffeomorphism 
$\varphi\col S \rightarrow S$ and
on the primitive
$N$--th root $q$ of
$1$. 
\end{thm}

Note that $L^q_\varphi$ is an
isomorphism of a vector space of
very large dimension
$N^{3g+p-3}$, and consequently
encodes a lot of information.
Extracting invariants from
$L^q_\varphi$ provides simpler
invariants of
$\varphi$, such as the
projectivized spectrum of
$L^q_\varphi$. We can also
normalize $L_\varphi$ so that it
has determinant 1,  in which case
its trace gives an invariant of
$\varphi$ defined up to
multiplication by a root of
unity. 

Explicit computations of these invariants in certain examples 
are provided in \cite{Liu2}. 

As is often the case with
invariants from Topological Quantum
Field Theory, the invariants
extracted from
$L^q_\varphi$ are by
themselves unlikely to have many
practical applications.  What is
more interesting is their
connections with other
combinatorial and geometric
objects. 

As this work was being developed, the type of functions occurring in explicit computations  hinted at a connection between the invariant of Theorem~\ref{thm:intro:Invariant}, the Kashaev $6j$--symbols developed in \cite{Kash1}, and the link invariants introduced by Kashaev \cite{Kash1, Kash2}, Baseilhac and Benedetti \cite{BasBen02, BasBen03, BasBen05, BasBen06}; see also Murakami-Murakami \cite{MM}. This connection has now been elucidated by the authors and Hua Bai \cite{Bai2, BBL}. Whereas the current article focuses on irreducible representations, \cite{BBL} investigates another type of representations of the quantum Teichm\"uller space, called local representations, which are somewhat simpler to analyze and more closely connected to the combinatorics of ideal triangulations. The classification of these local representations follows the same lines as the classification of irreducible representations, in terms of complex edge weights for ideal triangulations. An analogue of Theorem~\ref{thm:intro:Invariant} then associates to a homotopically aperiodic diffeomorphism $\varphi\col S \to S$ a large matrix $K_\varphi^q$, well-defined up to conjugation and multiplication by a root of unity. If one decomposes a local representation into its irreducible components, the  invariant $L_\varphi^q$ of Theorem~\ref{thm:intro:Invariant} and its generalizations discussed in Section~\ref{sect:Invariants} occur as building blocks of this $K_\varphi^q$. It can then be shown that the trace of $K_\varphi^q$ coincides with the invariant that, following the original insights of Kashaev,  Baseilhac and Benedetti \cite{BasBen05} associate to the hyperbolic metric of the mapping torus $M_\varphi$. A crucial step \cite{Bai2} is an explicit identification between the intertwining operator that a local representation associates to a diagonal exchange, and the $6j$--symbols that Kashaev defines using the representation theory of the Weyl Hopf algebra.

The results of this paper are very reminiscent of a well-known principle in quantum algebra, which is that the representations of a quantum group are in  correspondence with representations of the original non-quantum Lie group or algebra. It would also be conceptually helpful to  establish a
connection with the quantum group
constructions of
\cite{BFK, FG, FGL}, or with the skein theory of \cite{Prz, PrzSi, Tur}.

\medskip

It is a pleasure to thank Hua Bai, 
Leonid Chekhov and Bob Penner for
very helpful conversations. In
particular, this work originated
from lectures given by Leonid
Chekhov at U.S.C., and the reader
familiar with \cite{FC} will easily
recognize our debt to the last
paragraph of that paper. We are also
grateful to Bob
Guralnick, Chuck Lanski, Susan
Montgomery and Lance Small for
algebraic consulting, and to the referee for misprint hunting. 

This work was partially 
supported by the grant
DMS-0103511 from the National
Science Foundation.

\section{The Chekhov-Fock algebra}
\label{sect:Chekhov-Fock}

Let $S$ be an oriented punctured
surface of finite topological
type, obtained by removing a
finite set
$\{v_1, v_2, \dots, v_p \}$ from
the closed oriented surface
$\bar S$. Let $\lambda$ be an
ideal triangulation of $S$,
namely the intersection with $S$
of the 1--skeleton of a
triangulation of $\bar S$ whose
vertex set is equal to
$\{v_1, v_2, \dots, v_p \}$. In
other words, $\lambda$ consists
of finitely many disjoint simple
arcs $\lambda_1$,
$\lambda_2$, \dots, $\lambda_n$
going from puncture to puncture
and decomposing $S$ into finitely
many triangles with vertices at
infinity. Note that $n=-
3\chi(S)=6g +3p-6$, where
$\chi(S)$ is the Euler
characteristic of $S$, $g$ is the
genus of $\bar S$ and $p$ is the
number of punctures of $S$. In
particular, we will require that
$p\geq 3$ when $g=0$ to guarantee
the existence of such ideal
triangulations.

The complement $S-\lambda$ has
$2n$ spikes converging towards
the punctures, and each spike is
delimited by one $\lambda_i$ on
one side and one $\lambda_j$ on
the other side, with possibly
$i=j$. For
$i$,
$j\in
\{1, \dots, n\}$, let
$a_{ij}$ denote the number of
spikes of $S-\lambda$ which are
delimited on the left by
$\lambda_i$ and on the right by
$\lambda_j$ as one moves towards
the end of the spike, and set
\begin{equation*}
\sigma_{ij}= a_{ij}-a_{ji}.
\end{equation*} Note that
$\sigma_{ij}$ can only belong to
the set $\{-2, -1, 0, +1, +2\}$,
and that $\sigma_{ji}=
-\sigma_{ij}$. 

In the shear coordinates for
Teichm\"uller space associated to
the ideal triangulation
$\lambda$, the antisymmetric
bilinear form with matrix
$(\sigma_{ij})$ is closely related
to the Weil-Petersson closed
2--form on Teichm\"uller space
$\mathcal T(S)$. Compare 
\cite{Pen92, PaPen, Bon,
SozBon}, according to the type
of Teichm\"uller space
considered.

 The
\emph{Chekhov-Fock algebra}
associated to the ideal
triangulation $\lambda$ is the
algebra $\mathcal T_\lambda^q$
defined by the generators
$X_i^{\pm1}$, with
$i=1$, $2$, \dots, $n$, and by the
skew-commutativity relations
\begin{equation*} X_iX_j =
q^{2\sigma_{ij}} X_jX_i
\end{equation*} for every $i$,
$j$ (in addition to the relations
$X_i X_i^{-1} = X_i^{-1}X_i=1$). 

In particular, the Chekhov-Fock
algebra
$\mathcal T_\lambda^q$  is
an iterated skew-polynomial
algebra (see \cite{Coh}) as well
as a  special type of
multiparameter quantum torus (see
\cite[Chap.~I.2] {BroGoo}). What
is really important here is that
its algebraic structure is
 tied to the combinatorics of the
ideal triangulation $\lambda$ of
the surface $S$. 

We first analyze the algebraic
structure of
$\mathcal T_\lambda^q$.

\section{The structure of the
Weil-Petersson form}
\label{sect:StructureWP}

The skew-commutativity
coefficients
$\sigma_{ij}$ form an
antisymmetric matrix
$\Sigma$, which defines an
antisymmetric bilinear form
$\sigma\col  \mathbb Z^n \times
\mathbb Z^n
\rightarrow \mathbb Z$. The key
technical step to understanding
the algebraic structure of
$\mathcal T_\lambda^q$ is to
classify the bilinear form
$\sigma$ over the integers. Recall
that two bilinear forms on
$\mathbb Z^n$, with respective
matrices $\Sigma$ and $\Sigma'$,
are \emph{equivalent over $\mathbb
Z$} if there exists a base change
matrix $A\in\mathrm{GL}_n(\mathbb
Z)$ such that $\Sigma'=A\Sigma
A^{\mathrm t}$.

\begin{prop}
\label{prop:StructureWP} The
antisymmetric bilinear form
$\sigma\col  \mathbb Z^n \times
\mathbb Z^n
\rightarrow \mathbb Z$ is
equivalent over $\mathbb Z$ to the
block diagonal form  consisting of
$g$ blocks 
$\left(\begin{matrix} 0&-2\\2&0
\end{matrix}\right)$, $k$ blocks 
$\left(\begin{matrix} 0&-1\\1&0
\end{matrix}\right)$ and $p$
blocks $\bigl(0\bigr)$ on the
diagonal, where $g$ is the genus
of $S$,
$p$ is its number of punctures,
and
$k=2g+p-3$. 
\end{prop}
\begin{proof} Let $\Gamma\subset
S$ be the graph dual to the ideal
triangulation
$\lambda$. Note that every vertex
of $\Gamma$ is trivalent, and
that $\Gamma$ is a deformation
retract of $S$. 

The coordinates of the $\mathbb
Z^n$ considered above correspond
to the components of $\lambda$.
In a more intrinsic way, we
consequently have a natural
isomorphism between this $\mathbb
Z^n$ and the group
$\mathcal H(\lambda; \mathbb Z)$
of all assignments of integer
weights
 to the components of
$\lambda$ or, equivalently, to
the edges of $\Gamma$. In
particular, $\sigma$ is now an
antisymmetric bilinear form on
$\mathcal H(\lambda; \mathbb Z)$.

We first give a homological
interpretation of $\mathcal
H(\lambda; \mathbb Z)$ and 
$\sigma$, as is now somewhat
standard when analyzing the Thurston
intersection form on train tracks
(see for instance
\cite{Bon}). 

Let $\widehat\Gamma$ be the
oriented graph obtained from
$\Gamma$ by keeping the same
vertex set and by replacing each
edge of
$\Gamma$ by two oriented edges
which have the same end points as
the original edge, but which have
opposite orientations. In
particular, every vertex of
$\widehat\Gamma$ now has valence
6. There is a natural projection
$p\col \widehat\Gamma \rightarrow
\Gamma$ which is one-to-one on the
vertex set of $\widehat\Gamma$ and
two-to-one on the interior of the
edges of $\widehat\Gamma$. 

There is a unique way to thicken
$\widehat\Gamma$ to a surface
$\widehat S$ such that:
\begin{enumerate}
\item $\widehat S$ deformation
retracts to $\widehat\Gamma$;
\item as one goes around a vertex
$\widehat v$  of $\widehat\Gamma$
in
$\widehat S$, the orientations of
the edges of $\widehat\Gamma$
adjacent to
$\widehat v$ alternately point
towards and away from $\widehat
v$;
\item the natural projection
$p\col \widehat\Gamma \rightarrow
\Gamma$ extends to a $2$--fold
branched covering $\widehat S
\rightarrow S$, branched along the
vertex set of
$\widehat\Gamma$.
\end{enumerate} Indeed, the last
two conditions completely
determine the local model for the
inclusion of
$\widehat\Gamma $ in $\widehat S$
near the vertices of
$\widehat\Gamma$. 

Let $\tau\col \widehat S
\rightarrow
\widehat S$ be the covering
involution of the branched
covering $p\col \widehat S
\rightarrow  S$. Note
that $\tau$ respects $\widehat
\Gamma$, and reverses the
orientation of its edges. 
 
\begin{lem} 
\label{lem:Homology} There is a
natural identification between
$\mathcal H(\lambda; \mathbb Z)
\cong\mathbb Z^n$ and the subgroup
of 
$H_1(
\widehat S ) =H_1(
\widehat S; \mathbb Z)$ consisting
of those
$\widehat\alpha$ such that
$\tau_*(
\widehat\alpha ) =
-\widehat\alpha$.  
\end{lem}

\begin{proof} Every assignment
$\alpha
\in
\mathcal H(\lambda; \mathbb Z)$
of weights to the edges of
$\Gamma$ lifts to a
$\tau$--invariant edge weight
assignment
$\widehat\alpha$  for
$\widehat\Gamma$. Because the
edges of $\widehat \Gamma$ are
oriented, $\widehat\alpha$
actually defines a 1--chain on
$\widehat\Gamma$, whose boundary
is equal to 0 because each edge
$\widehat e$ of
$\widehat\Gamma$ is paired with
the edge $\tau(\widehat
e)$ which has the same
$\widehat\alpha$--weight but such
that $\partial  \tau(\widehat
e) = -\partial\widehat e$.
Therefore, we can interpret 
$\widehat\alpha $ as an element
of $ H_1( \widehat\Gamma
)$. Note that
$\tau_*(
\widehat\alpha )=
-\widehat\alpha$ since $\tau$
reverses the orientation of the
edges of $\widehat\Gamma$.

Conversely, every $\widehat\alpha
\in H_1( \widehat\Gamma
)$ associates an integer
weight to each edge of
$\widehat\Gamma$, by considering
its algebraic intersection number
with an arbitrary point in the
interior of the edge. If in
addition 
$\tau_*(
\widehat\alpha ) =
-\widehat\alpha$, this defines a
$\tau$--invariant edge weight
system on $\widehat\Gamma$, and
therefore an element of $\mathcal
H(\lambda; \mathbb Z)$.  

This identifies
$\mathcal H(\lambda; \mathbb Z)$
to the set of those
$\widehat\alpha \in H_1(
\widehat\Gamma )
= H_1( \widehat S
)$ such that
$\tau_*(
\widehat\alpha ) =
-\widehat\alpha$. 
\end{proof}

\begin{lem} 
\label{lem:Intersection} If
$\alpha$,
$\beta \in \mathcal H(\lambda;
\mathbb Z)$ correspond to
$\widehat\alpha$,
$\widehat\beta \in  H_1(
\widehat S )$  as
in Lemma~\ref
{lem:Homology}, then $\sigma(
\alpha, \beta ) $ is equal to the
algebraic intersection number
$\widehat\alpha \cdot
\widehat\beta$. 
\end{lem}

\begin{proof} It suffices to
check this for each generator
$\alpha_i\in \mathcal H(\lambda;
\mathbb Z)$ assigning weight 1 to
the edge
$e_i$ of $\Gamma$ dual to the
component
$\lambda_i$ of $\lambda$, and
weight 0 to the other edges of
$\Gamma$.  By definition,
$\sigma(\alpha_i,
\alpha_j)=\sigma_{ij}$ is equal
to the number of times $e_i$
appears to the immediate left
(as seen from the vertex)  of
$e_j$ at a vertex of $\Gamma$,
minus the number of times $e_i$
appears to the immediate right of
$e_j$. The corresponding homology
class $\widehat\alpha_i \in
H_1(
\widehat S )$ is
realized by the oriented closed
curve $c_i$ that is the union of
the two oriented edges of
$\widehat\Gamma$ lifting $e_i$.
In particular, $c_i$ and $c_j$
meet only at vertices of
$\widehat\Gamma$ corresponding to
common end points of the edges
$e_i$ and $e_j$ in $\Gamma$. When
$e_i$ is immediately to the left
of $e_j$ at a vertex of $\Gamma$,
it easily follows from our
requirement that edge
orientations alternately point in
and out at the vertices of
$\widehat\Gamma$ that the
corresponding intersection
between $c_i$ and $c_j$ has
positive sign. Similarly, an end
of $e_i$ which is immediately to
the right of an end of $e_j$
contributes a $-1$ to the
algebraic intersection number of
$c_i$ with $c_j$. It follows that
$\sigma(\alpha_i, \alpha_j) = c_i
\cdot c_j = \widehat\alpha_i
\cdot \widehat\alpha_j$. 

Therefore, $\sigma( \alpha, \beta
) =\widehat\alpha \cdot
\widehat\beta$ for every 
$\alpha$,
$\beta \in \mathcal H(\lambda;
\mathbb Z)$. 
\end{proof}

We now analyze in more detail the
branched covering $p\col 
\widehat S
\rightarrow S$. We claim that the
covering is trivial near the
punctures of $S$. Indeed, if
$\widehat C$ is a simple closed
curve going around a puncture in
$\widehat S$, the collapsing of
$\widehat S$ to $\widehat\Gamma$
sends $\widehat C$ to a curve
which is oriented by the
orientation of the edges of
$\Gamma$. This follows from our
requirement that the orientations
alternately point in and out at
each vertex of $\widehat\Gamma$.
Since the covering involution
$\tau$ reverses the orientation
of the edges of $\widehat\Gamma$,
we conclude that $\tau$ respects
no puncture of $\widehat S$. In
other words, a puncture of $S$
lifts to two distinct punctures
of $\widehat S$, and the covering
is trivial on a neighborhood of
this puncture.

The branched covering $p\col 
\widehat S
\rightarrow S$  is classified by a
homomorphism
$\pi_1(S-V)
\rightarrow \mathbb Z/2$, where
$V$ is the set of branch
points of $p$, namely the vertex
set of
$\Gamma$. Since the covering is
trivial near the punctures of
$S$, the  corresponding class
$H^1(S-V; \mathbb Z/2)$ is dual
to the intersection with $S-V$ of
a 1--submanifold $K\subset V$ with
$\partial K=V$. 
One can arrange by surgery that
$K$ consists only of arcs. Let 
$D\subset S$ be a disk containing
$K$, and let $\widehat D$ be its
preimage in $\widehat S$. The main
point here is that the restriction
$\widehat S -\widehat D
\rightarrow S-D$ is now a trivial
unbranched covering. In
particular, $\widehat S$ is the
union of $\widehat D$ and of two
copies $\widehat S_1$ and
$\widehat S_2$ of $S-D$. 

The restriction of $p$ to
$\widehat D\rightarrow D$ is a
2--fold branched covering of a
disk, with $4g +2p-4 $ branch
points. It follows that $\widehat
D$ is a surface of genus
$k=2g+p-3$ with two boundary
components. In addition, the covering involution is conjugate to a hyperelliptic involution of $\widehat D$, so that
the induced homomorphism $\tau_*$ acts on
$H_1(\widehat D;\mathbb Z)$  by
multiplication by $-1$.

Let $\widehat D^0$, $\widehat
S_1^0$ and $\widehat S_2^0$ be the
closed surfaces obtained by capping
off the punctures and boundary
components of $\widehat D$,
$\widehat S_1$ and $\widehat
S_2$, respectively. In addition,
for $i=1$, 2, \dots, $p$, let $C_i$
be a small curve going
counterclockwise around the $i$--th
puncture in $S$, and let
$\widehat C_{i1}$ and
$\widehat C_{i2}$ be its respective
lifts in $\widehat S_1$ and
$\widehat S_2$. Then
$H_1(\widehat S)$ is isomorphic to
$H_1(\widehat D^0) \oplus
H_1(\widehat S_1^0) \oplus
H_1(\widehat S_2^0) \oplus V$,
where $V$ is the subgroup generated 
the $\widehat C_{i1}$
and $\widehat C_{i2}$. Note that
the only relation between the
homology classes of these $2p$
curves is that they add up to 0, so
that
$V\cong \mathbb Z^{2p-1}$. 

Lemma~\ref{lem:Homology}
 identifies the space
$\mathcal H(\lambda; \mathbb Z)$
of edge weight assignments to the
subspace $\{
\widehat\alpha \in  H_1(\widehat
S) ;
\tau_*(\widehat\alpha)=
-\widehat\alpha\}$. By
construction, the
isomorphism $\tau_*$ of
$H_1(\widehat S) \cong H_1(\widehat 
D^0)
\oplus H_1(\widehat S_1^0) \oplus
H_1(\widehat S_2^0) \oplus V$ acts
by multiplication by
$-1$ on $H_1(\widehat 
D^0)$, exchanges the two factors 
$H_1(\widehat S_1^0) \cong
H_1(\widehat S_2^0) \cong
H_1(\bar S)$, and acts on $V$ by
transposing each pair $\{
\widehat C_{i1}, \widehat C_{i2}
\}$. (Recall that $\bar S$ is the
closed surface such that $S =
\bar S- \{v_1, v_2, \dots, v_p \}$.)
It follows that
$\mathcal H(\lambda; \mathbb Z)$
consists of those
$(x,y,-y,z)$ in $H_1(\widehat S)
\cong H_1(\widehat  D^0)
\oplus H_1( \bar S) \oplus
H_1( \bar S) \oplus V$ such that
$\tau_*(z)=-z$. This provides an
isomorphism $\mathcal H(\lambda;
\mathbb Z)
\cong H_1(\widehat D_0) 
\oplus H_1(\bar S)
\oplus W$, where $W=\{ z\in V;
\tau_*(z)=-z\} \cong \mathbb Z^p$.

By Lemma~\ref
{lem:Intersection}, the
bilinear form
$\sigma$ is the restriction to
$\mathcal H(\lambda; \mathbb Z)$
of the intersection form of
$H_1(\widehat S)$. We
conclude that the three factors
of the decomposition 
 $\mathcal H(\lambda; \mathbb Z)
\cong H_1(\widehat D_0) 
\oplus H_1(\bar S)
\oplus W$ are orthogonal for
$\sigma$, that the restriction of
$\sigma$ to $H_1(\widehat D_0)$ is
the intersection form of $\widehat
D_0$, that its restriction to 
$H_1(\bar S)$ is \emph{twice} the
intersection form of $\bar S$
(because $y\in H_1(\bar S)$ lifts
to $(0,y,-y,0)
\in H_1(\widehat S) \cong
H_1(\widehat D_0) 
\oplus H_1(\bar S)
\oplus H_1(\bar S)
\oplus V $), and that
$\sigma$ is 0 on
$W\cong \mathbb Z^p$. 

Since $\widehat D_0$ and
$\bar S$ are closed surfaces of
respective genus $k$ and $g$,
this concludes the proof of
Proposition~\ref{prop:StructureWP}.
\end{proof}

A consequence of
Proposition~\ref{prop:StructureWP} is
that the kernel of the bilinear
form $\sigma$, namely
\begin{equation*}
\mathrm{Ker}\, \sigma = \{ \alpha
\in \mathcal H(\lambda; \mathbb
Z);\, \forall \beta 
\in \mathcal H(\lambda; \mathbb
Z),\, \sigma (\alpha, \beta)=0 \},
\end{equation*}
 is isomorphic to $\mathbb Z^p$. We
can precise this result as follows. 
Index the punctures of
$S$ from 1 to $p$. For $i=1$,
\dots, $p$ and $j=1$, \dots, $n$,
let 
$k_{ij}
\in
\{0,1,2\}$ denote the number of
ends of the component $\lambda_j$
of $\lambda$ that converge to the
$i$--th puncture. Note that
$\sum_{i=1}^p (k_{i1}, k_{i2},
\dots, k_{in}) = (2,2,
\dots, 2)$ since each
$\lambda_j$ has two ends. 

\begin{lem}
\label{lem:StructureKerWP} In $
\mathcal H(\lambda; \mathbb
Z)\cong
\mathbb Z^n$, the kernel
$\mathrm{Ker}\, \sigma$ is the
abelian subgroup freely
generated by the $p$ vectors
$(1,1,\dots, 1)$ and
$(k_{i1}, k_{i2},
\dots, k_{in})$, for $i=1$, \dots
$p-1$. 
\end{lem}
\begin{proof}   Using the notation
of the proof of
Proposition~\ref{prop:StructureWP},
$\mathrm{Ker}\, \sigma$ corresponds
to the subspace
$W$ of $\mathcal H(\lambda;
\mathbb Z)
\cong H_1(\widehat D_0) 
\oplus H_1(\bar S)
\oplus W$. We
need to backtrack through the
definition of $W$.

Recall that, for each $i=1$, \dots,
$p$, we picked an oriented closed
curve
$C_i$ going counterclockwise around
the
$i$--th puncture of $S$, and that
we lifted it to curves
$\widehat C_{i1}$ and
$\widehat C_{i2}$ in
$\widehat S_1$ and $\widehat S_2$,
respectively. The only relation
between the $\widehat C_{i1}$ and
$\widehat C_{j2}$ is that their
sum is 0, so that they generate a subspace
$V\cong \mathbb Z^{2p-1}$ of
$H_1(\widehat S)$. Then, $W$
consists of those $z\in V$ such that
$\tau_*(z)=-z$. 

Since $\tau$ exchanges $\widehat C_{i1}$ and
$\widehat C_{i2}$, it follows
that $W$ is the abelian subgroup
freely generated by the $\widehat
C_{i1} -\widehat C_{i2}$, for
$i=1$, \dots, $p-1$, and by the
element
$H = \sum_{i=1}^p \widehat
C_{i1}= -\sum_{i=1}^p \widehat
C_{i2}$. 

As we retract the surface
$\widehat S$ to the graph
$\widehat \Gamma$, the curves
$\widehat C_{i1}$ and
$\widehat C_{i2}$ are sent to
curves in $\widehat \Gamma$
which, because of the alternating
condition for the edge
orientations at the vertices of
$\widehat \Gamma$, either follow
the orientation of the edges of
$\widehat\Gamma$ or go  against
this orientation everywhere. In
addition, because
$\tau$ reverses the orientation
of $\widehat\Gamma$, exactly one
of these two curves follow the
orientation. It follows that, for
the identifications
$\mathcal H(\lambda; \mathbb Z)
\cong \mathbb Z^n$ and  $\mathcal
H(\lambda; \mathbb Z)
\cong 
\{
\widehat\alpha \in  H_1(\widehat
S; \mathbb Z) ;\,
\tau_*(\widehat\alpha)=
-\widehat\alpha\}$, the vector 
$(k_{i1}, k_{i2},
\dots, k_{in})$ corresponds to
$\varepsilon_i (\widehat C_{i1}
-\widehat C_{i2}) \in W\subset 
H_1(\widehat S; \mathbb Z)$,
where $\varepsilon_i=+1$ when
$C_{i1}$ is sent to an
orientation preserving curve of
$\widehat\Gamma$, and
$\varepsilon_i=-1$ otherwise.
Note that what
determines $\varepsilon_i$ is our
choice of the disk $D \subset
\bar S$ in the proof of
Proposition~\ref{prop:StructureWP}.

Because each component
$\lambda_i$ of $\lambda$ has two
ends, 
\begin{equation*} (1,1,\dots, 1) =
\frac 12 \sum_{i=1}^p  (k_{i1},
k_{i2},
\dots, k_{in}),
\end{equation*} and it follows
that $(1,1, \dots, 1)\in \mathbb
Z^n$ corresponds to 
\begin{equation*}
\frac 12 \sum_{i=1}^p
\varepsilon_i (\widehat C_{i1}
-\widehat C_{i2})=
\varepsilon_p H
+\sum_{i=1}^{p-1} \delta_i  
(\widehat C_{i1} -\widehat C_{i2})
\end{equation*} with $\delta_i =
\frac{\varepsilon_i -
\varepsilon_p}2 = \pm1$. Since 
$\mathrm{Ker}\, \sigma =W $ is
freely generated by 
$H $ and by the
$\widehat C_{i1} -\widehat
C_{i2}$, for
$i=1$, \dots, $p-1$, it follows
that it is also generated by
those elements that, for the
identification $\mathcal
H(\lambda; \mathbb Z)
\cong \mathbb Z^n$, correspond to 
$(1,1,\dots, 1)$ and
$(k_{i1}, k_{i2},
\dots, k_{in})$, for $i=1$, \dots
$p-1$. 
\end{proof}

For a positive integer $N$, we
will also need to consider the
\emph{$N$--kernel} of $\sigma$,
defined as
\begin{equation*}
\mathrm{Ker}_N\, \sigma = \{
\alpha
\in \mathcal H(\lambda; \mathbb
Z); \,\forall \beta 
\in \mathcal H(\lambda; \mathbb
Z), \,\sigma (\alpha, \beta)
\in N\mathbb Z \}.
\end{equation*} Note that
$\mathrm{Ker}_N\,
\sigma$ contains $\mathcal
H(\lambda;N\mathbb Z)$. It
therefore makes sense to consider
its image in $ \mathcal H(\lambda;
\mathbb Z)/\mathcal H(\lambda;
N\mathbb Z)= \mathcal H(\lambda;
\mathbb Z_N)$, where $\mathbb
Z_N$ denotes the cyclic group
$\mathbb Z/ N\mathbb Z$.

\begin{lem}
\label{lem:StructureKerWPodd}
When $N$ is odd, the $N$--kernel 
$\mathrm{Ker}_N\,
\sigma$ is equal to the preimage
in  $\mathcal H(\lambda;\mathbb
Z)$ of the
$\mathbb Z_N$--submodule of 
$\mathcal H(\lambda;
\mathbb Z_N)\cong (\mathbb
Z_N)^n$ freely generated by  the
$p$ vectors
$(1,1,\dots, 1)$ and
$(k_{i1}, k_{i2},
\dots, k_{in})$, for $i=1$, \dots
$p-1$.
\end{lem}

\begin{proof} The image of the
$N$--kernel 
$\mathrm{Ker}_N\,
\sigma$ is the kernel
$\mathrm{Ker}\, \bar
\sigma$ of the form $\bar
\sigma\col  \mathcal H(\lambda; 
\mathbb Z_N) \times
\mathcal H(\lambda; \mathbb Z_N)
\rightarrow
\mathbb Z_N$ induced by
$\sigma$.  Replacing the
coefficient ring
$\mathbb Z$ by $\mathbb Z_N$, the
proof of
Proposition~\ref{prop:StructureWP}
provides an isomorphism
$\mathcal H(\lambda; \mathbb Z_N)
\cong H_1(\widehat D_0;
\mathbb Z_N) 
\oplus H_1(\bar S; \mathbb Z_N)
\oplus W_N$, where $W_N$ is the
image of the subspace $W$. The
three factors
$H_1(\widehat D_0;
\mathbb Z_N) $, 
$H_1(\bar S; \mathbb Z_N)$ and $
W_N$ are orthogonal for
$\bar\sigma$, and the restriction
of
$\bar\sigma$ to each factor is
the intersection form of
$\widehat D_0$, twice the
intersection form of $\bar S$,
and 0, respectively. 

Because $N$ is odd, $2$ is
invertible in $\mathbb Z_N$. If
follows that 
$\mathrm{Ker}\, \bar
\sigma = W_N$. The proof of
Lemma~\ref{lem:StructureKerWP}
now shows that
$W_N$ is freely generated by
$(1,1,\dots, 1)$ and by 
$(k_{i1}, k_{i2},
\dots, k_{in})$, for $i=1$, \dots
$b-1$.
\end{proof}

When $N$ is even,
$\mathrm{Ker}_N\,
\sigma$ contains additional
elements. Let $\alpha_1$,
$\alpha_2$, \dots, $\alpha_{2g}$
form a basis for $H_1(\bar S;
\mathbb Z _2)$. We can represent
$\alpha_i$ by a family $a_i$ of
curves immersed in the graph
$\Gamma \subset S$ dual to
$\lambda$ and passing at most
once across each edge of
$\Gamma$. Let
$l_{ij} \in
\{0,1\}$ be the number of times
$a_i$ traverses the $j$-th edge
of $\Gamma$.

\begin{lem}
\label{lem:StructureKerWPeven}
When $N$ is even, the $N$--kernel 
$\mathrm{Ker}_N\,
\sigma$ is equal to the preimage
in  $\mathcal H(\lambda;\mathbb
Z)$ of the direct sum
$A \oplus B \subset\mathcal
H(\lambda;
\mathbb Z_N)$ of the $\mathbb
Z_N$--submodule $A\cong
(\mathbb Z_N)^{g}$ freely generated
by  the vectors
$(1,1,\dots, 1)$ and
$(k_{i1}, k_{i2},
\dots, k_{in})$ for $i=1$, \dots
$p-1$, and of the
submodule $B \cong
(\mathbb Z_2)^{2g}$
generated by the $(l_{j1}\frac{N}2,
l_{j2}\frac{N}2,
\dots, l_{jn}\frac{N}2)$ with
$j=1$, \dots
$2g$.
\end{lem}

\begin{proof} 
The difference with
Lemma~\ref{lem:StructureKerWPeven}
is that 
$2\frac N2=0$ in $\mathbb Z_N$.
Therefore, in  $\mathcal H(\lambda;
\mathbb Z_N)
\cong H_1(\widehat D_0;
\mathbb Z_N) 
\oplus H_1(\bar S; \mathbb Z_N)
\oplus W_N$, the kernel
$\mathrm{Ker}\, \bar
\sigma$ is now the direct sum
$A\oplus B'$ of
$A=W_N$ and of the subspace $B'$ of
$H_1(\bar S;\mathbb Z_N)$
consisting of those elements which
are divisible by
$\frac N2$. As before
$A=W_N \cong
(\mathbb Z_N)^{g}$ is freely generated
by
$(1,1,\dots, 1)$ and by 
$(k_{i1}, k_{i2},
\dots, k_{in})$, for $i=1$, \dots
$p-1$. 

The factor $B'$ is also the image
$\frac N2 H_1(\bar S; \mathbb Z_2)$
of the group homomorphism $H_1(\bar
S; \mathbb Z_2) \to H_1(\bar S;
\mathbb Z_N)$ defined by
multiplication by $\frac N2$. To
identify explicit generators for
$B'\subset\mathcal
H(\lambda;
\mathbb Z_N)$, it is convenient to consider
the \emph{transfer map} $T\col 
H_1(S;\mathbb Z_2)
\rightarrow   H_1(\widehat
S;\mathbb Z_2)$, which to a cycle
in $S$ associates its preimage in
$\widehat S$. Its image is
contained in $\{ \alpha \in
H_1(\widehat S ; \mathbb Z_2) ;
\tau_*(\alpha)=\alpha
\}\cong \mathcal H(\lambda;\mathbb
Z_2)$. If $\alpha_i'\in H_1
(S;\mathbb Z_2)$ is represented by
the above family of curves $a_i$, it
is immediate from definitions that
$T(\alpha_i')$ corresponds to the
vector $(l_{j1}, l_{j2}, \dots,
l_{jn})$  in $\mathcal H(\lambda;\mathbb
Z_2) \cong (\mathbb Z_2)^n$. 

In the set-up of
Proposition~\ref{prop:StructureWP},
the transfer map $T$  can be
geometrically realized by
representing a class
$\alpha\in H_1(S;\mathbb Z_2)$ by
a curve $a$ contained in $S-D$;
then
$T(\alpha)$ is the class of $a_1
+a_2$, where $a_1$ and $a_2$ are
copies of $a$ in the two copies
$S_1$ and $S_2$  of $S-D$ contained
in
$\widehat S$. In particular, if we
start with a class $\alpha \in
H_1(\bar S; \mathbb Z_2)$, lift it
to a class $\alpha' \in H_1(
S; \mathbb Z_2)$ and consider
its image $T(\alpha') \in \mathcal
H(\lambda; \mathbb Z_2)
\cong H_1(\widehat D_0;
\mathbb Z_2) 
\oplus H_1(\bar S; \mathbb Z_2)
\oplus W_2$, the projection of
$T(\alpha')$ to the factor $H_1(\bar
S; \mathbb Z_2)$ is exactly equal
to $\bar\alpha$. As a consequence,
 $H_1(\bar S; \mathbb Z_2)
\oplus W_2$ is isomorphic to
$B_2
\oplus W_2$ if $B_2\subset\mathcal
H(\lambda;
\mathbb Z_2)$ denotes the subspace
generated by the $T(\alpha_i')$.

Multiplying everything by $\frac
N2$ we conclude that, in 
$\mathcal H(\lambda; \mathbb Z_N)
\cong H_1(\widehat D_0;
\mathbb Z_N) 
\oplus H_1(\bar S; \mathbb Z_N)
\oplus W_N$, the kernel
$\mathrm{Ker}\, \bar\sigma
=0\oplus \frac N2 H_1(\bar S;
\mathbb Z_2)
\oplus W_N$ is equal to
$0\oplus B\oplus W_N$ where $B=\frac
N2 B_2$ is generated by the 
vectors 
$(l_{j1}\frac{N}2,
l_{j2}\frac{N}2,
\dots, l_{jn}\frac{N}2)$.
\end{proof}

\section{The algebraic structure
of the Chekhov-Fock algebra}
\label{sect:StructureCF}

\begin{lem}
\label{lem:basisCF} The monomials
$X_1^{k_1} X_2^{k_2}
\dots X_n^{k_n}$, with $k_1$,
$k_2$, \dots, $k_n \in \mathbb Z$,
form a basis for $\mathcal
T_\lambda^q$, considered as a
vector space. 
\end{lem}
\begin{proof} This immediately
follows from the fact that 
$\mathcal T_\lambda^q$ is an
iterated (Laurent) skew-polynomial
algebra, and can also be described
as the vector space freely
generated by these monomials and
endowed with the appropriate
multiplication. See
\cite[\S 2.1]{Coh} or \cite[\S
1.7]{Kas}.
\end{proof}

\begin{thm}
\label{prop:structureCF} The
Chekhov-Fock algebra $\mathcal
T_\lambda^q$ is isomorphic to the
algebra $\mathcal W_{g,k,p}^q$
defined by generators
$U_i^{\pm1}$, $V_i^{\pm1}$, with
$i=1$, \dots, $g+k$, and
$Z_j^{\pm1}$ with $j=1$, \dots,
$p$ and by the following
relations:
\begin{enumerate}
\item each $U_i$ commutes with all
generators except $V_i^{\pm1}$;
\item each $V_i$ commutes with all
generators except $U_i^{\pm1}$;
\item $U_iV_i=q^4 V_iU_i$ for
every $i=1$, \dots, $g$;
\item $U_iV_i=q^2 V_iU_i$ for
every $i=g+1$, \dots, $g+k$;
\item each $Z_j$ commutes with
all generators.
\end{enumerate} Here $g$ is the
genus of the surface $S$, $p$ is
its number of punctures and
$k=2g+p-3$. In addition, the
isomorphism between
$\mathcal T_\lambda^q$ and
$\mathcal W_{g,k,p}^q$ can be
chosen to send monomial to
monomial. 
\end{thm}

\begin{proof} Let $F_n$ be the
free group generated by the set
$\{X_1,
\dots, X_n \}$. We can rephrase
the definition of $\mathcal
T_\lambda^q$ by saying that it is
is the quotient of the group
algebra $\mathbb C[F_n]$ by the
2--sided ideal generated by all
elements $ X_iX_j -
q^{\sigma_{ij}}X_jX_i$. 

Note that the abelianization of
$F_n$ is canonically isomorphic to
$\mathbb Z^n$. In addition, if we
identify two words $a$, $b\in
F_n$ to their images in $\mathcal
T_\lambda^q$ and if $\bar a$ and
$\bar b$ denote their images in
$\mathbb Z^n$, then
$ba=q^{\sigma(\bar a,\bar b)}ab$
in
$\mathcal T_\lambda^q$. 

Consider the base change isomorphism
$\mathbb Z^n
\rightarrow \mathbb Z^n$  provided
by Proposition~\ref
{prop:StructureWP}, under which
$\sigma$ becomes block diagonal.
Lift this isomorphism to a group
isomorphism
$ F_n
\rightarrow F_n$, which itself
induces an algebra isomorphism\break
$\Phi\col  \mathbb C [F_n]
\rightarrow \mathbb C[F_n]$. If
we denote the generators of the
first $F_n$ by\break
 $\{U_1, V_1, U_2,
V_2, \dots, U_{g+k}, V_{g+k},
Z_1,  Z_2, \dots, Z_p\}$,
it immediately follows from
definitions that
$\Phi$ induces an isomorphism
from $\mathcal W_{g,k,p}^q$ to
$\mathcal T_\lambda^q$. This
isomorphism sends monomial to
monomial since it comes from an
isomorphism of
$F_n$.
\end{proof}

The monomials
$aX_1^{i_1}X_2^{i_2} \dots
X_n^{i_n}$, with $i_j \in \mathbb
Z$ and  $a \in
\mathbb C$, play a particularly
important r\^ole in the structure
of $\mathcal T^q_\lambda$ and of
its representations. Let
$\mathcal M^q_\lambda$ denote the
set of all such monomials that are
different from 0. The
multiplication law of
$\mathcal T^q_\lambda$ induces a
group law on $\mathcal
M^q_\lambda$. 

The elements $aX_1^0X_2^0
\dots X_n^0$ form a subgroup of
$\mathcal M^q_\lambda$ isomorphic
to the multiplicative group
$\mathbb C^* = \mathbb C-\{0\}$.
There is also a natural group
homomorphism 
$\mathcal M^q_\lambda \rightarrow
\mathbb Z^n = \mathcal
H(\lambda;\mathbb Z)$ which to 
$X=a X_1^{i_1}X_2^{i_2} \dots
X_n^{i_n}$ associates the vector  
$\bar X = (i_1, i_2,
\dots, i_n)$. This defines a
central extension 
\begin{equation*} 1 \rightarrow 
\mathbb C^* \rightarrow
\mathcal M^q_\lambda \rightarrow
\mathbb Z^n \rightarrow 1
\end{equation*} whose algebraic
structure is completely
determined by the commutation
property that
$XY= q^{2\sigma(\bar X,
\bar Y)} YX$ for every $X$,
$Y\in
\mathcal M^q_\lambda$.

Let $\mathcal Z^q_\lambda$ be the
center of $\mathcal M^q_\lambda$.
An immediate consequence of
Lemma~\ref{lem:basisCF} is that
the center of the algebra
$\mathcal T^q_\lambda$ consists of
all sums of elements of 
$\mathcal Z^q_\lambda$.  We now
analyze the structure of
$\mathcal Z^q_\lambda$.

We first introduce preferred
elements of $\mathcal
Z^q_\lambda$. By
Lemma~\ref{lem:StructureKerWP},
$\mathcal Z^q_\lambda$ contains
the element $X_1X_2 \dots X_n$.
However, it is better to introduce
its scalar multiple  
\begin{equation*} H= q^{
-\sum_{i<i'}
\sigma_{ii'}} X_1X_2
\dots X_n.
\end{equation*}

Similarly,
Lemma~\ref{lem:StructureKerWP}
shows that the center
$\mathcal Z^q_\lambda$ contains
the element\break
$X_1^{k_{i1}} X_2^{k_{i2} }\dots
X_n^{k_{in}}\in \mathcal
T_\lambda^q$ associated to the
$i$--th puncture of $S$, where
$k_{ij}
\in
\{0,1,2\}$ denotes the number of
ends of the component $\lambda_j$
of $\lambda$ that converge to this
$i$--th puncture.   Again, we
consider 
\begin{equation*} P_i= q^{
-\sum_{j<j'} k_{ij} k_{ij'}
\sigma_{jj'}} 
 X_1^{k_{i1}} X_2^{k_{i2} }\dots
X_n^{k_{in}}
\end{equation*}

The  $q$--factors in the 
definition of
$H$ and of the
$P_i$ are specially defined to
guarantee invariance under
re-indexing of the $X_j$. This
choice of scalar factors is
classically known as the Weyl
quantum ordering.

\begin{lem} 
\label{lem:PHrelations} For every
integer $N$, 
\begin{align*} H^2 &=P_1P_2 \dots
P_p \\ H^N &= q^{ -N^2\sum_{i<i'}
\sigma_{ii'}} X_1^N X_2^N \dots
X_n^N \\  P_i^N &= q^{
-N^2\sum_{j<j'} k_{ij} k_{ij'}
\sigma_{jj'}} 
 X_1^{Nk_{i1}} X_2^{Nk_{i2} }\dots
X_n^{Nk_{in}}
\end{align*}
\end{lem} 
\begin{proof}  The $P_i$ and $H$
belong to the subset $\mathcal
A\subset
\mathcal Z_\lambda^q$ consisting
of all elements  of the form
\begin{equation*} q^{-\sum_{j<k}
\sigma_{i_ji_{k}}} X_{i_1}X_{i_2}
\dots X_{i_m}
\end{equation*}  Note that the
fact that the elements of
$\mathcal A$ are central implies
that
$\sum_k
\sigma_{ji_k}=0$ for every $j$.
It immediately follows that, for
every
$A$ and $B\in \mathcal A$, the
product $AB$ is also in $\mathcal
A$. Also, an element of
$\mathcal A$ is invariant under
permutation of the $X_{i_j}$ (and
subsequent adjustment of the
$q$--factor). 

The three equations of
Lemma~\ref{lem:PHrelations}
immediately follow from these
observations, using for the first
equation the fact that
$\sum_i k_{ij}=2$ for every $j$. 
\end{proof}

\begin{prop}
\label{prop:centerCF} When $q$ is
not a root of unity, the center
$\mathcal Z^q_\lambda$ of the
monomial group $\mathcal
M^q_\lambda$ is equal to the
direct sum of $\mathbb C^*$ and
of the  abelian
 subgroup freely generated (as an
abelian group) by the above
elements $H$ and
$P_i$ with $i=1$, \dots, $p-1$.
\end{prop}
\begin{proof} This immediately
follows from the algebraic
structure of $\mathcal
M^q_\lambda$ and from
Lemma~\ref{lem:StructureKerWP}. 
\end{proof}

When $q^2$ is a primitive $N$--th
root of unity, the center  
$\mathcal Z^q_\lambda$ contains
additional elements, such as the
$X_i^{ N}$. 
Lemma~\ref{lem:PHrelations}
provides relations between $H^N$,
the $X_i^N$ and the $P_j^N$.

\begin{prop}
\label{prop:centerCFodd} If $q^2$
is a primitive
$N$--th root of unity with $N$
odd, the center $\mathcal
Z^q_\lambda$ of the monomial group
$\mathcal M^q_\lambda$ is
generated by the $X_i^{ N}$ with
$i=1$,
\dots, $n$, by the element
$H$, and by the $P_j$ with
$j=1$, \dots, $p-1$.

 In addition, if $W$ denotes the
direct sum of 
$\mathbb C^*$ and of the free
abelian group generated by the
$X_i^{ N}$,
$H$ and $P_j$, with $i=1$,
\dots, $n$ and
$j=1$, \dots,
$p-1$, then $\mathcal
Z^q_\lambda$ is
isomorphic to the quotient of $W$ by
the relations:
\begin{align*}
 H^N&= q^{ -N^2\sum_{i<i'}
\sigma_{ii'}} X_1^NX_2^N
\dots X_n^N\\
 P_j^N&=q^{ -N^2\sum_{k<k'}
k_{jk} k_{jk'}
\sigma_{kk'}} (X_1^N)^{k_{j1}}
(X_2^N)^{k_{j2}}
\dots (X_n^N)^{k_{jn}}.
\end{align*}
\end{prop}

\begin{proof} Again, this
immediately follows from our
analysis of $\mathrm{Ker}\,
\sigma_N$ in Lemma~\ref
{lem:StructureKerWPodd}, together
with the relations of Lemma~\ref
{lem:PHrelations}.
\end{proof}

It should be noted that, when
$q^2$ is an $N$--th root of
unity, then $q^N=\pm1$ so that
the $q$--factors in the relations
of Proposition~\ref
{prop:centerCFodd} are  equal to
$\pm1$. In later sections, we
will choose
$q$ so that these factors are
actually equal to 1, making these
relations less intimidating.

When $N$ is even, the structure
of $\mathrm{Ker}\,
\sigma_N$ is more complicated, and
consequently so is the structure
of  $\mathcal Z^q_\lambda$. Let
$\alpha_1$,
$\alpha_2$, \dots, $\alpha_{2g}$
form a basis for $H_1(\bar S;
\mathbb Z _2)$. We can represent
$\alpha_k$ by a family $a_k$ of
curves immersed in the graph
$\Gamma \subset S$ dual to
$\lambda$ and passing at most
once across each edge of
$\Gamma$. Let
$l_{ki} \in
\{0,1\}$ be the number of times
$a_k$ traverses the $i$-th edge
of $\Gamma$. Define 
\begin{equation*} A_k =
q^{-\frac{N^2}4\sum_{i<i'}
l_{ki}l_{ki'}\sigma_{ii'}}
X_1^{\frac N2 l_{k1}} X_2^{\frac
N2 l_{k2}}
\dots X_n^{\frac N2 l_{kn}} \in
\mathcal T^q_\lambda. 
\end{equation*}

As in
Lemma~\ref{lem:PHrelations}, 
\begin{align*} A_k^2 &=
q^{-{N^2}\sum_{i<i'}
l_{ki}l_{ki'}\sigma_{ii'}} X_1^{N
l_{k1}} X_2^{N l_{k2}}
\dots X_n^{N l_{kn}}\\  &= X_1^{N
l_{k1}} X_2^{N l_{k2}}
\dots X_n^{N l_{kn}}
\end{align*} since
$q^{N^2}=(\pm1)^N=1$ because
$N$ is even.

\begin{prop}
\label{prop:centerCFeven} If
$q^2$ is a primitive
$N$--th root of unity with $N$
even, the center $\mathcal
Z^q_\lambda$ of the monomial
group $\mathcal M^q_\lambda$ is 
generated by
$\mathbb C^*$, by the $X_i^{ N}$
with
$i=1$,
\dots, $n$, by the element
$H$, by the $P_j$ with
$j=1$, \dots, $p-1$, and by the
$A_k$ with
$k=1$, \dots, $p-1$.

  In addition, if $W$ denotes the
direct sum of $\mathbb
C^*$ and of the free abelian group
generated by the $X_i^{N}$,
$H$, $P_j$ and
$A_k$, with
$i=1$,
\dots, $n$, $j=1$, \dots,
$p-1$ and
$k=1$, \dots,
$2g$, then $\mathcal
Z^q_\lambda$ is 
isomorphic to the quotient of $W$ by
the relations:
\begin{align*}
 H^N &=  X_1^NX_2^N
\dots X_n^N\\
 P_j^N&= (X_1^N)^{k_{j1}}
(X_2^N)^{k_{j2}}
\dots (X_n^N)^{k_{jn}}\\ A_k^2&= 
(X_1^N)^{l_{k1}} (X_2^N)^{l_{k2}}
\dots (X_n^N)^{l_{kn}}
\end{align*}

\end{prop}
\begin{proof} Again, this follows
from Lemma~\ref
{lem:StructureKerWPeven},
together with the relations of
Lemma~\ref {lem:PHrelations} and
the fact that $q^{N^2}=1$ when
$N$ is even.
\end{proof}

\section{Finite-dimensional
representations of the
Chekhov-Fock algebra}
\label{sect:RepTheory}

This section is devoted to the
classification of the
finite-dimensional representations
of the algebra $\mathcal
T^q_\lambda$, namely of the
algebra homomorphisms 
$\rho\col \mathcal T^q_\lambda
\rightarrow
\mathrm{End} (V)$ from
$\mathcal T^q_\lambda$ to the
algebra of endomorphisms of a
finite-dimensional vector space
$V$ over $\mathbb C$. Recall that
two such representations 
$\rho\col \mathcal T^q_\lambda
\rightarrow
\mathrm{End} (V)$ and 
$\rho'\col \mathcal T^q_\lambda
\rightarrow
\mathrm{End} (V')$ are
\emph{isomorphic} if there exists
a linear isomorphism
$L\col V\rightarrow V'$ such that
$\rho'(X) = L \cdot \rho(X) \cdot
L^{-1}$ for every $X \in \mathcal
T^q_\lambda$, where $\ \cdot\ $
denotes the composition of maps
$V' \rightarrow V \rightarrow V
\rightarrow  V'$.  Also, 
$\rho\col
\mathcal T^q_\lambda
\rightarrow
\mathrm{End} (V)$ is
\emph{irreducible} if it does not
respect any proper subspace
$W\subset V$.

Having determined the algebraic
structure of
$\mathcal T^q_\lambda$ in
Section~\ref{sect:StructureCF},
the classification of its
representations is an easy
exercise (see
Lemmas~\ref{lem:RepWeyl},
\ref{lem:RepProduct} and
\ref{lem:RepLaurent}). The main
challenge is to state this
classification in an intrinsic
way which is tied to the topology
of the ideal triangulation
$\lambda$. This is done in
Theorem~\ref{thm:ClassRepCF} in a
first step, and then in
Theorems~\ref{thm:ClassRepCFodd}
and \ref{thm:ClassRepCFeven} in a
more concrete way. 

It is not hard to see that the
Chekhov-Fock algebra $\mathcal
T^q_\lambda$ cannot admit any
finite-dimensional representation
unless $q$ is a root of unity. In
this case, our results will
heavily depend on the number $N$
such that
$q^2$ is a primitive
$N$--th root of unity. 

In addition to the structure
theorems of
Section~\ref{sect:StructureCF},
our analysis of the
representations of $\mathcal
T^q_\lambda$ is based on the
following elementary (and
classical) facts. 

\begin{lem}
\label{lem:RepWeyl} Let $\mathcal
W^q$ be the algebra defined by
the generators
$U^{\pm1}$, $V^{\pm1}$ and by the
relation $UV=q^2VU$. If $q^2$ is
a primitive $N$--th root of
unity, every irreducible
representation of $\mathcal W^q$
has dimension $N$, and is
isomorphic to a representation
$\rho_{uv}$ defined by 
\begin{equation*}
\rho_{uv}(U)=u\left(
\begin{matrix} 1 & 0 & 0  &
\dots & 0& 0\\ 0 & q^2 & 0
&\dots  &0& 0\\ 0 & 0 & q^4
&\dots & 0& 0\\
\dots & \dots & \dots &\dots &
\dots& \dots\\  
\dots & \dots & \dots &\dots &
\dots& \dots\\ 0 & 0 & 0 &\dots  &
q^{2N-4}&0\\ 0 & 0 & 0 & \dots &
0 & q^{2N-2}
\end{matrix}
\right)
\end{equation*} and 
 \begin{equation*}
\rho_{uv}(V)=v\left(
\begin{matrix} 0 & 0 & 0 & \dots&
0 & 1\\ 1 & 0 & 0 &\dots & 0 & 0\\
0 & 1 & 0 & \dots& 0 & 0\\
\dots & \dots & \dots & \dots&
\dots & \dots\\
 \dots & \dots & \dots & \dots&
\dots & \dots\\ 0 & 0 & 0 &
\dots& 0 & 0\\ 0 & 0 & 0 &\dots &
1 & 0
\end{matrix}
\right)
\end{equation*} for some $u$,
$v\in \mathbb C-\{0\}$. In
addition, two such
representations $\rho_{uv}$ and
$\rho_{u'v'}$ are isomorphic if
and only if $u^N= (u')^N$ and
$v^N=(v')^N$. 
\end{lem}
\begin{proof} Note that $U^N$ and
$V^N$ are central in
$\mathcal W^q$. If
$\rho$ is an irreducible
representation, it must
consequently send $U^N$ to a
homothety $u_1\, \mathrm{Id}$ and
$V^N$ to a homothety $v_1\,
\mathrm{Id}$. In addition,
$\rho(V)$  sends an eigenvector
of $\rho(U)$ corresponding to an
eigenvalue
$\ell$ to another eigenvector of
$\rho(U)$ corresponding to the
eigenvalue $\ell q^2$. It easily
follows that $\rho$ is isomorphic
to a representation
$\rho_{uv}$ for some $u$, $v$
such that $u^N=u_1$ and
$v^N=v_1$. 

If the representations
$\rho_{uv}$ and $\rho_{u'v'}$ are
isomorphic, then necessarily $u^N=
(u')^N$ and
$v^N=(v')^N$ by consideration of
the homotheties $\rho_{uv}(U^N)$,
$\rho_{u'v'}(U^N)$,
$\rho_{uv}(V^N)$ and
$\rho_{u'v'}(V^N)$.  Conversely,
conjugating
$\rho_{uv}$ by the isomorphism
$\rho_{uv}(U)$ gives the
representation 
$\rho_{u'v'}$ with $u'=u$ and
$v'=vq^2$; it follows that the
isomorphism class of
$\rho_{uv}$ depends only on $u$
and $v^N$. Similarly, the
representation obtained by
conjugating
$\rho_{uv}$ by the isomorphism
$\rho_{uv}(V)$ is equal to the
representation 
$\rho_{u'v'}$ with $u'=uq^2$ and
$v'=v$. It follows that the
isomorphism class of
$\rho_{uv}$ depends only on $u^N$
and $v^N$.
\end{proof}

\begin{lem}
\label{lem:RepProduct} Let $q^2$
be a primitive $N$--th root of
unity, and let $\mathcal W^q$ be
the algebra defined by the
generators
$U^{\pm1}$, $V^{\pm1}$ and by the
relation $UV=q^2VU$. Let
$\mathcal W$ be any algebra. Any
irreducible finite-dimensional
representation of the tensor
product $\mathcal W \otimes
\mathcal W^q$ is isomorphic to
the tensor product
$\rho_1 \otimes \rho_2\col 
\mathcal W \otimes \mathcal W^q
\rightarrow \mathrm{End} (W_1
\otimes W_2)$ of two irreducible
representations $\rho_1\col 
\mathcal W \rightarrow
\mathrm{End}(W_1)$ and
$\rho_2\col 
\mathcal W^q \rightarrow
\mathrm{End}(W_2)$. Conversely,
the tensor product of two such
irreducible representations is
irreducible. 
\end{lem}

\begin{proof} Consider an
irreducible representation
$\rho\col 
\mathcal W \otimes \mathcal W^q
\rightarrow \mathrm{End} (W)$,
with $W$ a finite-dimensional
vector space over $\mathbb C$.
Let $W_1 \subset W$ be an
eigenspace of $\rho(1\otimes U)$,
corresponding to the eigenvalue
$u$. Then $\rho(1 \otimes V^i)$
sends $W_1$ to the eigenspace
$W_{i+1}$ of $\rho(1\otimes U)$
corresponding to the eigenvalue
$uq^{2i}$. Also, $\mathcal
W\otimes 1$ commutes with
$1\otimes U$, and 
$\rho(\mathcal W\otimes 1)$
consequently preserves each
$W_i$. Noting that
$\rho(1\otimes V^N)$ is a
homothety since $1\otimes V^N$ is
central, it follows that
$\bigoplus_{i=1}^N W_i$ is
invariant under $\rho (\mathcal W
\otimes \mathcal W^q)$ , and is
therefore equal to $W$ by
irreducibility of $\rho$. 

If $\rho(\mathcal W\otimes 1)$
respected a proper subspace
$W_1'$ of $W_1$, then by the
above remarks the subspace 
$\bigoplus_{i=1}^N
\rho(1\otimes V^i)(W_1')$ would be
a proper subspace invariant under
$\rho (\mathcal W \otimes
\mathcal W^q)$. By irreducibility
of $\rho$, it follows that the
representation 
$\rho_1\col 
\mathcal W \rightarrow
\mathrm{End}(W_1)$ defined by
restriction of $\rho (\mathcal W
\otimes 1)$ to $W_1$ is
irreducible. 

All the pieces are now here to
conclude that the representation
$\rho$ of $\mathcal W \otimes
\mathcal W^q$ over
$W=\bigoplus_{i=1}^N W_i$ is
isomorphic to the tensor  product
of $\rho_1\col 
\mathcal W \rightarrow
\mathrm{End}(W_1)$ and of a
representation 
$\rho_2\col 
\mathcal W^q \rightarrow
\mathrm{End}(W_2)$ of the type
described in
Lemma~\ref{lem:RepWeyl}. 

Conversely,  consider the tensor
product
$\rho$ of two  irreducible
representations
$\rho_1\col 
\mathcal W \rightarrow
\mathrm{End}(W_1)$ and
$\rho_2\col 
\mathcal W^q \rightarrow
\mathrm{End}(W_2)$, where
$\rho_2$ is as in 
Lemma~\ref{lem:RepWeyl}. Let
$L_u\subset W_2$ be the
(1--dimensional) eigenspace of
$\rho_2(U)$ corresponding to the
eigenvalue
$u$, so that $W_1 \otimes L_u$ is
the eigenspace of $\rho(1\otimes
U)$ corresponding to the
eigenvalue
$u$. If $W'
\subset W_1\otimes W_2$ is
invariant under $\rho$, in
particular it is invariant under
$\rho (1 \otimes \mathcal W^q)$,
and it follows from
Lemma~\ref{lem:RepWeyl} that
$W' \cap (W_1 \otimes L_u)$ is
non-trivial since
$\rho(1\otimes U^N) = u^N
\mathrm{Id}$. The subspace 
$W' \cap (W_1 \otimes L_u)$ is
also invariant under
$\rho(\mathcal W\otimes 1)$, and
must therefore be equal to all of
$W_1 \otimes L_u$ by
irreducibility of $\rho_1$. 
Therefore, $W'$ contains $W_1
\otimes L_u$, from which it
easily follows that
$W'=W_1 \otimes W_2$. This proves
that $\rho$ is irreducible. 
\end{proof}

\begin{lem}
\label{lem:RepLaurent} Let 
$\mathbb C[Z^{\pm1}]$ be   the
algebra of Laurent polynomials in
the variable $Z$, and let
$\mathcal W$ be any algebra. Any
irreducible finite-dimensional
representation of the tensor
product $\mathcal W \otimes
\mathbb C[Z^{\pm1}]$ is
isomorphic to the tensor product
$\rho_1 \otimes \rho_2\col 
\mathcal W \otimes 
\mathbb C[Z^{\pm1}]
\rightarrow \mathrm{End} (V_1
\otimes V_2)$ of two irreducible
representations $\rho_1\col 
\mathcal W \rightarrow
\mathrm{End}(W_1)$ and
$\rho_2\col 
\mathbb C[Z^{\pm1}] \rightarrow
\mathrm{End}(W_2)$. Conversely,
the tensor product of two such
irreducible representations is
irreducible. 
\end{lem}
\begin{proof} This immediately
follows from the fact that $Z$ is
central in $\mathcal W \otimes
\mathbb C[Z^{\pm1}]$, and from
the fact that every irreducible
representation $\rho_2\col 
\mathbb C[Z^{\pm1}] \rightarrow
\mathrm{End}(W_2)$ has dimension
1 and is classified by the number
$z\in \mathbb C^*$ such that
$\rho_2(Z)= z\,
\mathrm{Id}_{W_2}$. 
\end{proof}

Recall that $\mathcal Z^q_\lambda$
denotes the center of the group
$\mathcal M^q_\lambda$ of non-zero
monomials in the Chekhov-Fock
algebra
$\mathcal T^q_\lambda$. 

Let $\rho\col \mathcal
T^q_\lambda \rightarrow
\mathrm{End} (V)$ be a
finite-dimensional irreducible
representation of
$\mathcal T^q_\lambda$. Every $X
\in \mathcal Z^q_\lambda$ is
central in
$\mathcal T^q_\lambda$, and its
image $\rho(X)$ consequently is a
homothety, namely of the form
$a\,
\mathrm{Id}_V$ for $a \in \mathbb
C$. We can therefore interpret
the restriction of $\rho$ to
$\mathcal Z^q_\lambda
\subset\mathcal T^q_\lambda$ as a
group homomorphism $\rho\col 
\mathcal Z^q_\lambda \rightarrow
\mathbb C^*$.  Note that
$\rho\col  
\mathcal Z^q_\lambda \rightarrow
\mathbb C^*$ coincides with the
identity on $\mathbb C^* \subset 
\mathcal Z^q_\lambda$.

\begin{thm}
\label{thm:ClassRepCF}
 Suppose that
$q^2$ is a primitive $N$--th root
of unity. Every irreducible
finite-dimensional 
representation
$\rho\col \mathcal T^q_\lambda
\rightarrow
\mathrm{End} (V)$ has dimension
$N^{3g+p-3}$ if $N$ is odd, and
$N^{3g+p-3}/2^g$ if $N$ is even
(where $g$ is the genus of the
surface $S$ and $p$ is its number
of punctures). Up to isomorphism,
$\rho$ is completely determined
by its restriction $\rho\col 
\mathcal Z^q_\lambda \rightarrow
\mathbb C^*$ to the center
$\mathcal Z^q_\lambda$ of the
monomial group $\mathcal
M^q_\lambda$ of $\mathcal
T^q_\lambda$. 

Conversely, every group
homomorphism $\rho\col  \mathcal
Z^q_\lambda \rightarrow
\mathbb C^*$ coinciding with the
identity on 
$\mathbb C^* \subset 
\mathcal Z^q_\lambda$ can be
extended to an irreducible
finite-dimensional  representation
$\rho\col \mathcal T^q_\lambda
\rightarrow
\mathrm{End} (V)$.
\end{thm}

\begin{proof} By Theorem~\ref
{prop:structureCF} and for
$k=2g+p-3$, the Chekhov-Fock
algebra
 $\mathcal T_\lambda^q$ is
isomorphic to the algebra
$\mathcal W_{g,k,p}^q$ defined by
generators
$U_i^{\pm1}$, $V_i^{\pm1}$, with
$i=1$, \dots, $g+k$, and
$Z_j^{\pm1}$ with $j=1$, \dots,
$p$ and by the following
relations:
\begin{enumerate}
\item each $U_i$ commutes with all
generators except $V_i^{\pm1}$;
\item each $V_i$ commutes with all
generators except $U_i^{\pm1}$;
\item $U_iV_i=q^4 V_iU_i$ for
every $i=1$, \dots, $g$;
\item $U_iV_i=q^2 V_iU_i$ for
every $i=g+1$, \dots, $g+k$;
\item each $Z_j$ commutes with
all generators.
\end{enumerate}  In particular,
$\mathcal T_\lambda^q$ is
isomorphic to the tensor product
of $g$ copies of the algebra
$\mathcal W^{q^2}$ (defined by
the generators
$U^{\pm 1}$,
$V^{\pm 1}$ and by the relation
$UV=q^4VU$), $k$ copies of the
algebra $\mathcal W^q$, and $p$
copies of the algebra $\mathbb
C[Z^{\pm 1}]$.  In addition, the
isomorphism 
$\mathcal W_{g,k,p}^q\cong
\mathcal T_\lambda^q$ can be
chosen to send the monomial group
of
$\mathcal W_{g,k,p}^q$ to the
monomial group $\mathcal
M_\lambda^q$  of
$\mathcal T_\lambda^q$.

By Lemmas~\ref{lem:RepProduct} and
\ref{lem:RepLaurent}, an
irreducible finite-dimensional
representation is therefore
isomorphic to a tensor product
$\rho_1 \otimes \rho_2 \otimes
\dots \otimes
\rho_{g+k+p}$ of irreducible
representations $\rho_i$ such
that $\rho_i$ is a representation
of $\mathcal W^{q^2}$ for $1\leq
i \leq g$, a representation of
$\mathcal W^{q}$ if $g+1 \leq i
\leq g+k$, and a representation
of $\mathbb C[Z]$ if $g+k+1\leq i
\leq g+k+p$. In particular, for
$g+k+1\leq i \leq g+k+p$, the
irreducible representation
$\rho_i$ must have dimension $1$,
and is determined by the complex
number $\rho(Z_i) \in \mathbb
C^*$. 

If $N$ is odd, then $q^2$ and
$q^4$ are both primitive $N$--th
roots of unity. It follows from
Lemma~\ref{lem:RepWeyl} that, for
$1\leq i \leq g+k$, the
representation $\rho_i$ has
dimension $N$ and is completely
determined by the two homotheties
$\rho_i(U_i^N)$ and
$\rho_i(V_i^N)$. As a consequence
$\rho$ has dimension
$N^{g+k}=N^{3g+p-3}$, as
announced, and is completely
determined by the homotheties
that are the images of
$U_i^N$,
$V_j^N$ and
$Z_l$. Since $U_i^N$,
$V_j^N$ and
$Z_l$ belong to the center of the
monomial group of
$\mathcal W_{g,k,p}^q\cong
\mathcal T_\lambda^q$, this shows
that $\rho$ is determined by the
restriction of
$\rho$ to this center $\mathcal
Z_\lambda^q$. 

When $N$ is even, then $q^2$ is a
primitive $N$--th root of unity,
but $q^4$ is a primitive $\frac
N2$--th root of unity. Lemma~\ref
{lem:RepWeyl} now implies that
$\rho_i$ has dimension $\frac N2$
if $i=1$, $2$, \dots, $g$, and
has dimension $N$ if $g+1\leq i
\leq g+k$. It follows that $\rho$
has dimension $(\frac N2)^gN^k=
N^{3g+p-3}/2^g$, as announced. In
addition, $\rho_i$ is determined
by the homotheties
$\rho\bigl (U_i^{\frac N2} \bigr)$
and
$\rho\bigl (V_i^{\frac N2}
\bigr)$ if
$i=1$,
$2$, \dots, $g$, and by
$\rho(U_i^N)$ and
$\rho(V_i^N)$ if $g+1\leq i \leq
g+k$. Consequently, $\rho$ is
completely determined by the
images of the
$U_i^{\frac N2}$,
$V_i^{\frac N2}$ with $1\leq i
\leq g$, of the $U_i^N$ and
$V_i^N$ with $g+1\leq i \leq
g+k$, and of the
$Z_i$ with $g+k+1\leq i \leq
g+k+p$. Since these elements all
belong to the center of the
monomial group of
$\mathcal W_{g,k,p}^q\cong
\mathcal T_\lambda^q$, this shows
that $\rho$ is determined by the
restriction of
$\rho$ to this center $\mathcal
Z_\lambda^q$. 

This concludes the proof of the
first statement of Theorem~\ref
{thm:ClassRepCF}. 

We prove the second statement
when $N$ is even. The odd case is
similar. 

Consider a group homomorphism 
$\rho\col  \mathcal Z^q_\lambda
\rightarrow
\mathbb C^*$ coinciding with the
identity on $\mathbb C^*$.
Lemma~\ref{lem:RepWeyl}
associates an irreducible
representation
$\rho_i$ of $\mathcal W^{q^2}$ to
the numbers 
$\rho(U_i^{\frac N2})$ and
$\rho(V_i^{\frac N2})$ when $1\leq
i \leq g$, an irreducible
representation
$\rho_i$ of $\mathcal W^{q}$ to
 $\rho(U_i^N)$ and $\rho(V_i^N)$
when $g+1\leq i \leq g+k$. When
$g+k+1\leq i \leq g+k+p$, there
is a 1--dimensional representation
$\rho_i$ of $\mathbb C
[Z^{\pm1}]$ such that
$\rho_i(Z)=\rho(Z_i)$. This
defines a representation
$\rho'=\rho_1 \otimes \rho_2
\otimes
\dots \otimes
\rho_{g+k+p}$ of 
$\mathcal W_{g,k,p}^q\cong
\mathcal T_\lambda^q$, which is
irreducible by
Lemma~\ref{lem:RepProduct}. It
remains to show that the group
homomorphism
$\rho'\col 
\mathcal Z^q_\lambda \rightarrow
\mathbb C^*$ induced by $\rho'$
coincides with the original group
homomorphism $\rho\col  \mathcal
Z^q_\lambda \rightarrow
\mathbb C^*$. But this
immediately follows from the fact
that the center of the monomial
group of $\mathcal
W_{g,k,p}^q\cong
\mathcal T_\lambda^q$ is the
product of $\mathbb C^*$ and of
the free abelian group generated
by the
$U_i^{\frac N2}$,
$V_i^{\frac N2}$ with $1\leq i
\leq g$, by the $U_i^N$ and
$V_i^N$ with $g+1\leq i \leq
g+k$, and by the
$Z_i$ with $g+k+1\leq i \leq
g+k+p$.

This concludes the proof, when
$N$ is even, of the property that
every $\rho\col  \mathcal
Z^q_\lambda \rightarrow
\mathbb C^*$ coinciding with the
identity on $\mathbb C^*$ can be
extended to an irreducible
representation $\rho=\rho'$ of
$\mathcal T_\lambda^q$. As
indicated above, the case where
$N$ is odd is almost identical. 
\end{proof}

To express Theorem~\ref
{thm:ClassRepCF} in a more
concrete and geometric way, we
now combine this result with
our analysis of the algebraic
structure of the center
$\mathcal Z^q_\lambda$ in
Propositions~\ref
{prop:centerCFodd} and
\ref{prop:centerCFeven}.

Recall that we associated the
element 
\begin{equation*} P_i = q^{
-\sum_{j<j'} k_{ij}
k_{ij'}\sigma_{jj'}} X_1^{k_{i1}}
X_2^{k_{i2} }\dots
X_n^{k_{in}}\in \mathcal
T_\lambda^q
\end{equation*}
 to the
$i$-th puncture of $S$, where
$k_{ij}
\in
\{0,1,2\}$ is the number of ends
of the component $\lambda_j$ of
$\lambda$ that converge to this
$i$--th puncture.  We also
considered the element
\begin{equation*} H= q^{
-\sum_{i<i'} \sigma_{ii'}}
X_1X_2\dots X_n.
\end{equation*}

\begin{thm}
\label{thm:ClassRepCFodd} If
$q^2$ is a primitive $N$--th root
of unity with $N$ odd, the
irreducible finite-dimensional
representation
$\rho\col \mathcal T^q_\lambda
\rightarrow
\mathrm{End} (V)$ is, up to
isomorphism, completely
determined by:
\begin{enumerate}
\item for $i=1$, $2$, \dots, $n$,
the number
$x_i \in \mathbb C^*$ such that
$\rho(X_i^N) = x_i\,
\mathrm{Id}_V$;
\item for $j=1$, $2$, \dots,
$p-1$, the $N$--th root $p_j$ of
$\varepsilon_j x_1^{k_{j1}}
x_2^{k_{j2}}
\dots x_n^{k_{jn}}$ such that 
$\rho(P_j) = p_j\,
\mathrm{Id}_V$;
\item the $N$--th root $h$ of
$\varepsilon_0 x_1 x_2 \dots x_n$
such that $\rho(H) = h\,
\mathrm{Id}_V$;
\end{enumerate} where
$\varepsilon_j= q^{ -N^2
\sum_{l<l'} k_{jl}
k_{jl'}\sigma_{ll'}}=\pm 1$ and 
$\varepsilon_0= q^{ -N^2
\sum_{l<l'} \sigma_{ll'}}=\pm 1$. 

Conversely, every such data of
numbers $x_i$, $p_j$ and $h\in
\mathbb C^*$ with 
$p_j^N=\break  \varepsilon_j x_1^{k_{j1}}
x_2^{k_{j2}}
\dots x_n^{k_{jn}}$ and
$h^N= \varepsilon_0 x_1 x_2 \dots
x_n$ can be realized by an
irreducible finite-dimensional
representation
$\rho\col \mathcal T^q_\lambda
\rightarrow
\mathrm{End} (V)$. 
\end{thm}

\begin{proof} Combine Theorem~\ref
{thm:ClassRepCF} and
Proposition~\ref
{prop:centerCFodd}.
\end{proof}

In the case where $N$ is even, we
had to use a basis
$\alpha_1$,
$\alpha_2$, \dots, $\alpha_{2g}$
for $H_1(\bar S; \mathbb Z_2)$.
After representing each
$\alpha_k$ by a family $a_k$ of
curves immersed in the graph
$\Gamma \subset S$ dual to
$\lambda$ and passing 
$l_{ki}
\in
\{0,1\}$ times across the
$i$--th edge of $\Gamma$, we
introduced the monomial 
\begin{equation*} A_k =
q^{-\frac{N^2}4\sum_{i<i'}
l_{ki}l_{ki'}\sigma_{ii'}}
X_1^{\frac N2 l_{k1}} X_2^{\frac
N2 l_{k2}}
\dots X_n^{\frac N2 l_{kn}} \in
\mathcal T^q_\lambda. 
\end{equation*}

\begin{thm}
\label{thm:ClassRepCFeven} If
$q^2$ is a primitive $N$--th root
of unity with $N$ even, the
irreducible finite-dimensional
representation
$\rho\col \mathcal T^q_\lambda
\rightarrow
\mathrm{End} (V)$ is, up to
isomorphism, completely
determined by:
\begin{enumerate}
\item for $i=1$, $2$, \dots, $n$,
the number
$x_i \in \mathbb C^*$ such that
$\rho(X_i^N) = x_i\,
\mathrm{Id}_V$;
\item for $j=1$, $2$, \dots,
$p-1$, the $N$--th root $p_j$ of
$x_1^{k_{j1}} x_2^{k_{j2}}
\dots x_n^{k_{jn}}$ such that 
$\rho(P_j) = p_j\,
\mathrm{Id}_V$;
\item the $N$--th root $h$ of
$x_1 x_2 \dots x_n$ such that
$\rho(H) = h\,
\mathrm{Id}_V$;
\item for $k=1$, $2$, \dots,
$2g$, the square root $a_k$ of 
$ x_1^{l_{k1}} x_2^{l_{k2}}
\dots x_n^{l_{kn}}$ such that
$\rho(A_k) = a_k\,
\mathrm{Id}_V$.
\end{enumerate}

Conversely, every such data of
numbers $x_i$, $p_j$, $h$ and
$a_k\in
\mathbb C^*$ with 
$p_j^N=  x_1^{k_{j1}} x_2^{k_{j2}}
\dots x_n^{k_{jn}}$,
$h^N= x_1 x_2 \dots x_n$ and
$a_k^2 = x_1^{l_{k1}} x_2^{l_{k2}}
\dots x_n^{l_{kn}}$ can be
realized by an irreducible
finite-dimensional representation
$\rho\col \mathcal T^q_\lambda
\rightarrow
\mathrm{End} (V)$. 
\end{thm}

\begin{proof} Combine Theorem~\ref
{thm:ClassRepCF} and
Proposition~\ref
{prop:centerCFeven}.
\end{proof}

\section{The quantum
Teichm\"uller space}
\label{sect:QuantumTeich}

As one moves from one ideal
triangulation  $\lambda$ of the
surface $S$ to another ideal
triangulation $\lambda'$, there
is a canonical isomorphism 
$\Phi_{\lambda\lambda'}^q \col 
\widehat {\mathcal T}_{\lambda'}^q
\rightarrow \widehat {\mathcal
T}_\lambda^q$ between the
fraction algebras of the
Chekhov-Fock algebras respectively
associated to two ideal
triangulations $\lambda$ and
$\lambda'$.

Here the \emph{fraction algebra} 
$\widehat {\mathcal T}_\lambda^q$
is the division algebra consisting
of all the formal fractions
$PQ^{-1}$ with
$P$, $Q \in {\mathcal
T}_\lambda^q$ and $Q\not =0$,
subject to the `obvious'
manipulation rules. In other
words,
$\widehat {\mathcal T}_\lambda^q$
is the division algebra of all
the non-commutative rational
fractions in the variables $X_i$,
subject to the relations $X_iX_j =
q^{2\sigma_{ij}}X_jX_i$. The
existence of such a fraction
algebra is guaranteed by the fact
that $\mathcal T_\lambda^q-
\{0\}$ satisfies the so-called
Ore condition in $\mathcal
T_\lambda^q$; see for instance
\cite{Coh, Kas}. 

The isomorphism
$\Phi_{\lambda\lambda'}^q \col 
\widehat {\mathcal T}_{\lambda'}^q
\rightarrow \widehat {\mathcal
T}_\lambda^q$ was introduced by
Chekhov and Fock \cite{FC} as a
quantum deformation of the
corresponding change of coordinates
in Thurston's shear coordinates for
Teichm\"uller space. See \cite{Liu}
for a version which is more detailed
(in particular with respect to
non-embedded diagonal exchanges)
and is better adapted to the
context of the current paper.  

To describe the isomorphism
$\Phi_{\lambda\lambda'}^q$, we
need to be a little more careful
with definitions. We will
henceforth agree that the data of
an ideal triangulation $\lambda$
also includes an indexing  of the
components
$\lambda_1$, $\lambda_2$, \dots,
$\lambda_n$ of $\lambda$ by the
set $\{ 1, 2, \dots, n\}$. Let
$\Lambda(S)$ denote the set of
isotopy classes of all such
(indexed) ideal triangulations of
$S$. 

The set $\Lambda(S)$ admits two
natural operations. The first one
is the \emph{re-indexing action}
of the permutation group
$\mathfrak S_n$, which to
$\lambda \in
\Lambda(S)$ and $\alpha \in
\mathfrak S_n$ associates the
indexed ideal triangulation
$\alpha\lambda$ whose
$i$--th component is equal to
$\lambda_{\alpha(i)}$. 

\begin{figure}[h]
\SetLabels
(.17 * .60) $\lambda_i$\\
(.17 * 1) $\lambda_j$\\
(.34 * .5) $\lambda_k$\\
(.17 * -.08) $\lambda_l$\\
(-.02 * .5) $\lambda_m$\\
(.84 * .60) $\lambda_i'$\\
(.84 * 1.02) $\lambda_j'$\\
(1.02 * .5) $\lambda_k'$\\
(.84 * -.08) $\lambda_l'$\\
(.65 * .5) $\lambda_m'$\\
\endSetLabels
\centerline{
\AffixLabels{\includegraphics{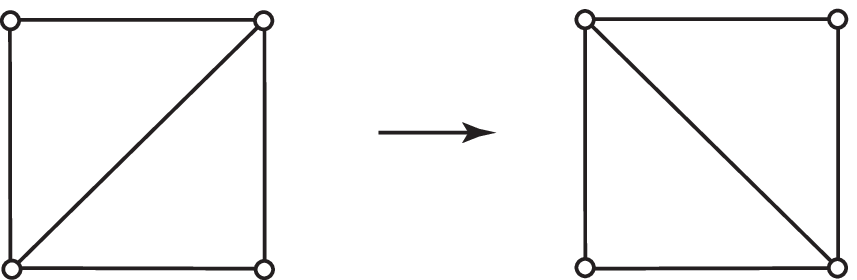}}
}
\caption{}
\label{pict:DiagEx}
\end{figure}

The second operation is the
\emph{$i$--th diagonal exchange}
$\Delta_i\col \Lambda(S)
\rightarrow
\Lambda(S)$ defined as follows.
In general, the $i$--th component
$\lambda_i$ of the ideal
triangulation
$\lambda$ separates two triangle
components $T_1$ and $T_2$ of
$S-\lambda$. The union $T_1 \cup
T_2 \cup \lambda_i$ is an open
square $Q$ with diagonal
$\lambda_i$. Then the ideal
triangulation $\Delta_i(\lambda)
\in \Lambda(S)$ is obtained from
$\lambda$ by replacing
$\lambda_i$ by the other diagonal
of the square $Q$, as in
Figure~\ref{pict:DiagEx}. This
operation is not defined when the
two sides of
$\lambda_i$ are in the same
component of $S-\lambda$, which
occurs when $\lambda_i$ is the
only component of $\lambda$
converging to a certain puncture;
in this case, we decide that
$\Delta_i(\lambda)=\lambda$. 

It may very well happen that two
distinct sides of the square $Q$
correspond to the same component
$\lambda_j$ of
$\lambda$. If, as in
Figure~\ref{pict:DiagEx}, we list
the components of
$\lambda$ in the boundary of $Q$
counterclockwise as
$\lambda_j$, $\lambda_k$,
$\lambda_l$ and $\lambda_m$, in
such a way that the diagonal
$\lambda_i$ goes from the
$\lambda_j\lambda_k$ corner to
the $\lambda_l\lambda_m$ corner,
we can list all possibilities as
follows, up to symmetries of the
square:
\begin{enumerate}
\item[1.\phantom{'}] The four
sides
$\lambda_j$, $\lambda_k$,
$\lambda_l$ and $\lambda_m$ of
the square $Q$ are all distinct;
in this case, we will say that
the diagonal exchange is
\emph{embedded}. 
\item[2.\phantom{'}]
$\lambda_j=\lambda_l$ and
$\lambda_k \not= \lambda_m$.
\item[2'.] $\lambda_k=\lambda_m$
and
$\lambda_j \not= \lambda_l$; note
that a diagonal exchange of this
type is the inverse of a diagonal
exchange of type 2.
\item[3.\phantom{'}]
$\lambda_j=\lambda_k$ and
$\lambda_l \not= \lambda_m$.
\item[3'.] $\lambda_j=\lambda_m$
and
$\lambda_k \not= \lambda_l$; note
that a diagonal exchange of this
type is the inverse of a diagonal
exchange of type 3.
\item[4.\phantom{'}]
$\lambda_j=\lambda_l$ and
$\lambda_k = \lambda_m$; note
that $S$ is the once punctured
torus in this case.
\item[5.\phantom{'}]
$\lambda_j=\lambda_k$ and
$\lambda_l = \lambda_m$; note
that $S$ is the three-times
punctured sphere in this case.
\item[5'.] $\lambda_j=\lambda_m$
and
$\lambda_k = \lambda_l$; note
that a diagonal exchange of this
type is the inverse of a diagonal
exchange of type 5.
\end{enumerate} Observe that
these different situations affect
the structure of $\mathcal
T_\lambda^q$ and 
$\mathcal T_{\lambda'}^q$ if
$\lambda'=\Delta_i(\lambda)$. For
instance, in $\mathcal
T_\lambda^q$,
$X_iX_j$ is equal to
$q^2 X_jX_i$ in Cases~1 and 2', is
equal to
$q^4 X_jX_i$ in Cases~2 and 4,
and is equal to $X_jX_i$ in
Cases~3, 3', 5 and 5'. Similarly,
in
$\mathcal T_{\lambda'}^q$,
$X_iX_j$ is equal to
$q^{-2} X_jX_i$ in Cases~1 and 2',
is equal to
$q^{-4} X_jX_i$ in Cases~2 and 4,
and is equal to $X_jX_i$ in
Cases~3, 3', 5 and 5'.

\begin{thm}[{\cite{FC, Liu}}]
\label{thm:CoordChange} There is
a unique family of isomorphisms 
$\Phi_{\lambda\lambda'}^q \col 
\widehat {\mathcal T}_{\lambda'}^q
\rightarrow \widehat {\mathcal
T}_\lambda^q$, indexed by pairs
of ideal triangulations
$\lambda$, $\lambda'\in
\Lambda(S)$, such that:
\begin{enumerate}
\item [(a)] for any 
$\lambda$,
$\lambda'$ and $\lambda'' \in
\Lambda(S)$,
$\Phi_{\lambda
\lambda''}^q = \Phi^q_{\lambda
\lambda'} \circ \Phi^q_{\lambda'
\lambda''};$
\item [(b)] if $\lambda'=\alpha
\lambda$ is obtained by
re-indexing $\lambda \in
\Lambda(S)$ by the permutation
$\alpha \in \mathfrak S_n$,
$\Phi_{\lambda\lambda'}^q$ is
defined by the property that
$\Phi_{\lambda\lambda'}^q(X_i) =
X_{\alpha(i)}$;
\item [(c)] if
$\lambda'=\Delta_i(\lambda)$ is
obtained from $\lambda$ by an
$i$--th diagonal exchange and if
we list all possible
configurations as in Cases~1-5'
above, then 
$\Phi_{\lambda\lambda'}^q$ is
defined by the property that 
$\Phi_{\lambda\lambda'}^q(X_h) =
X_h$ for every $h\not\in
\{i,j,k,l, m\}$, 
$\Phi_{\lambda\lambda'}^q(X_i) =
X_i^{-1}$, and:
\begin{enumerate}

\item  [(i)] in Case~1,
$$\begin{aligned}
\Phi_{\lambda\lambda'}^q(X_j) &=
(1+qX_i)X_j,\quad
&\Phi_{\lambda\lambda'}^q(X_k) &=
(1+qX_i^{-1})^{-1} X_k, \\
\Phi_{\lambda\lambda'}^q(X_l) &=
(1+qX_i)X_l, 
&\Phi_{\lambda\lambda'}^q(X_m) &=
(1+qX_i^{-1})^{-1} X_m;
\end{aligned}$$ 

\item [(ii)] in Case~2,
$$\begin{aligned}
\Phi_{\lambda\lambda'}^q(X_j) &=
(1+qX_i)(1+q^3X_i)X_j,\\
\Phi_{\lambda\lambda'}^q(X_k) &=
(1+qX_i^{-1})^{-1} X_k, \quad 
 \Phi_{\lambda\lambda'}^q(X_m) =
(1+qX_i^{-1})^{-1} X_m;
\end{aligned}$$ 

\item [(iii)] in Case~3,
$$\begin{aligned}
\Phi_{\lambda\lambda'}^q(X_j) &=
X_iX_j,
\quad\Phi_{\lambda\lambda'}^q(X_l)
&= (1+qX_i)X_l ,\\ 
\Phi_{\lambda\lambda'}^q(X_m) &=
(1+qX_i^{-1})^{-1} X_m;
\end{aligned}$$ 
 
\item [(iv)] in Case~4,
$$\begin{aligned}
\Phi_{\lambda\lambda'}^q(X_j) &=
(1+qX_i)(1+q^3X_i)X_j, \\
\Phi_{\lambda\lambda'}^q(X_k) &=
(1+qX_i^{-1})^{-1}
(1+q^3X_i^{-1})^{-1}X_k;
\end{aligned}$$

\item [(v)] in Case~5, 
$$\Phi_{\lambda\lambda'}^q(X_j) =
X_iX_j,\qquad
\Phi_{\lambda\lambda'}^q(X_l) =
X_iX_l.\eqno{\qed}$$ 
\end{enumerate}
\end{enumerate}
\end{thm}

The uniqueness of
$\Phi_{\lambda\lambda'}^q$ in
Theorem~\ref{thm:CoordChange}
immediately comes from the fact
that  any two ideal
triangulations $\lambda$ and
$\lambda'$ of $S$ can be
connected by a finite sequence of
diagonal moves and re-indexings
(see for instance \cite{Pen} for
this property). The difficult part
is to show that the isomorphism
$\Phi_{\lambda\lambda'}^q$  so
defined is independent of the
choice of the sequence of
diagonal moves and re-indexings.

The isomorphisms 
$\Phi_{\lambda\lambda'}^q \col 
\widehat {\mathcal T}_{\lambda'}^q
\rightarrow \widehat {\mathcal
T}_\lambda^q$  enable us to
associate an algebraic object to
the surface
$S$  in a way which does not
depend on the choice of an ideal
triangulation $\lambda$. For
this, consider the set of all
pairs $(X,\lambda)$ where
$\lambda \in \Lambda(S)$ is an
ideal triangulation of $S$ and
where
$X\in \widehat{\mathcal
T}^q_\lambda$. Define the
\emph{quantum Teichm\"uller
space}, as
\begin{equation*}
\widehat{\mathcal T}^q_S =
\{  (X,\lambda); \lambda \in
\Lambda(S), X\in \widehat{\mathcal
T}^q_\lambda
\}/\sim
\end{equation*} where the
equivalence relation
$\sim$ identifies 
 $(X, \lambda)$ to
$(X', \lambda')$ when
$X=\Phi^q_{\lambda\lambda'}(X')$.
 The set
$\widehat{\mathcal T}^S_\lambda$
inherits a natural division
algebra structure from that of
the $\widehat{\mathcal
T}^q_\lambda$. In fact, for any
ideal triangulation $\lambda$,
there is a natural isomorphism
between $\widehat{\mathcal
T}^q_S$ and $\widehat{\mathcal
T}^q_\lambda$.

The terminology is motivated by
the non-quantum (also called
semi-classical) case where
$q=1$ (see 
\cite{FC, Liu}, and compare
Section~\ref{sect:HyperShadow}). 
Consider the enhanced
Teichm\"uller space
$\mathcal T(S)$ of
$S$,  where each element consists
of a complete hyperbolic metric
defined up to isotopy together
with an orientation for each end
of $S$ that has infinite area for
the metric. Thurston's shear
coordinates for Teichm\"uller
space (see for instance
\cite{Bon, FC, Liu},
\cite{Thu86} for a dual version, and \cite{FoG} for a generalization)
associate to the  ideal
triangulation
$\lambda\in
\Lambda(S)$ a diffeomorphism
$\varphi_\lambda \col\mathcal T(S)
\rightarrow \mathbb R^n$. The
corresponding coordinate changes
$\varphi_{\lambda'}
\circ\varphi_\lambda^{-1}
\col \mathbb R^n \rightarrow
\mathbb R^n$ are rational
functions and, for the natural
identifications $\widehat{\mathcal
T}^1_\lambda\cong
\widehat{\mathcal T}^1_{\lambda'}
\cong \mathbb C(X_1, X_2, \dots,
X_n)$, it turns out that the
isomorphism $\mathbb C(X_1, X_2, 
\dots, X_n) \rightarrow \mathbb
C(X_1, X_2, \dots, X_n)$ induced
by $\varphi_{\lambda'}
\circ\varphi_\lambda^{-1}$
exactly coincides with 
$\Phi_{\lambda\lambda'}^1 \col 
\widehat {\mathcal T}_{\lambda'}^1
\rightarrow \widehat {\mathcal
T}_\lambda^1$. As a consequence,
there is a natural notion of
rational functions on $\mathcal
T(S)$, and the algebra of these
rational functions is naturally
isomorphic to 
$\widehat{\mathcal T}^1_S$. 

For a general $q$, the division
algebra
$\widehat{\mathcal T}^q_S$  can
therefore be considered as a
deformation of the algebra
$\widehat{\mathcal T}^1_S$ of all
rational functions on the
enhanced Teichm\"uller
space $\mathcal T(S)$. See
\cite{ Liu}. 

By analogy with the non-quantum
situation, we can think of the
natural isomorphism
$\widehat{\mathcal T}^p_\lambda
\rightarrow \widehat{\mathcal
T}^p_S$ as a parametrization of 
$\widehat{\mathcal T}^q_S$ by the
explicit algebra
$\widehat{\mathcal T}^q_\lambda$
associated to the ideal
triangulation $\lambda$. Pursuing
the analogy, we will call the
isomorphism 
$\Phi_{\lambda\lambda'}^q \col 
\widehat {\mathcal T}_{\lambda'}^q
\rightarrow \widehat {\mathcal
T}_\lambda^q$ the
\emph{coordinate change
isomorphisms} associated to the
ideal triangulations $\lambda$ and
$\lambda'$. 

Hua Bai \cite{Bai} proved that
the formulas of Theorem~\ref
{thm:CoordChange} are essentially
the only ones for which the
property holds, once we require
the $\Phi_{\lambda\lambda'}^q$ to
satisfy a small number of natural
conditions. In particular, the
quantum Teichm\"uller space is
a combinatorial object naturally
associated to the 2--skeleton of
the Harer-Penner simplicial
complex \cite{Har, Pen} of ideal
cell decompositions of $S$.

For future reference, we note:

\begin{lem} [\cite{Liu}]
\label{lem:CoordChange&HP} For
any two ideal triangulations
$\lambda$, $\lambda'$, the
coordinate change isomorphism 
$\Phi_{\lambda\lambda'}^q \col 
\widehat {\mathcal T}_{\lambda'}^q
\rightarrow \widehat {\mathcal
T}_\lambda^q$ sends the central
elements $H$, $P_1$, $P_2$,
\dots, $P_p$ of $\widehat
{\mathcal T}_{\lambda'}^q$ to 
the central elements $H$, $P_1$,
$P_2$,
\dots, $P_p$ of $\widehat
{\mathcal T}_{\lambda}^q$,
respectively. \qed
\end{lem}

As a consequence, $H$ and the
$P_i$ give well-defined central
elements of the quantum
Teichm\"uller space
$\widehat{\mathcal T}^q_S$, as
well as of its polynomial core
${\mathcal T}^q_S$ defined in
the next section.

\section{The polynomial core of
the quantum Teichm\"uller space}
\label{sect:PolCore}

The division algebras 
$\widehat{\mathcal T}^q_\lambda$
and $\widehat{\mathcal T}^q_S$
have a major drawback. They do
not admit any finite-dimensional
representations. Indeed, if there
was such a finite-dimensional
representation $\rho\col
\widehat {\mathcal T}^q_\lambda
\rightarrow
\mathrm{End}(V)$, then
$\rho(Q)\in \mathrm{End}(V)$
would be invertible for every $Q
\in \widehat {\mathcal
T}^q_\lambda -\{0\}$, by
consideration of
$\rho(Q^{-1})$. However,  since
${\mathcal T}^q_\lambda$ is
infinite-dimensional and
$\mathrm{End}(V)$ is
finite-dimensional, the
restriction 
$\rho\col {\mathcal T}^q_\lambda
\rightarrow
\mathrm{End}(V)$ has a huge
kernel, which provides many $Q$
for which $\rho(Q)=0$ is
non-invertible.

On the other hand, we saw in
\S\ref{sect:RepTheory} that the
Chekhov-Fock algebra
$\mathcal T_\lambda^q$ admits a
rich representation theory. This
leads us to introduce the
following definition. 

Let the
\emph{polynomial core} 
${\mathcal T}^q_S$ of the quantum
Teichm\"uller space
$\widehat{\mathcal T}^q_S$ be the
family $\{ \mathcal T_\lambda^q\}
_{\lambda \in \Lambda(S)}$ of all
Chekhov-Fock algebras $\mathcal
T_\lambda^q$, considered as
subalgebras of $\widehat{\mathcal
T}^q_S$, as $\lambda$ ranges over
all ideal triangulations of the
surface $S$. 

Given two ideal triangulations
$\lambda$ and $\lambda'$ and two
finite-dimensional
representations 
$\rho_\lambda\col \mathcal
T^q_\lambda
\rightarrow
\mathrm{End} (V)$ and 
$\rho_{\lambda'}\col \mathcal
T^q_{\lambda'}
\rightarrow
\mathrm{End} (V)$ of the
associated Chekhov-Fock algebras,
we would like to say that the two
representations correspond to
each other under the coordinate
change isomorphism 
$\Phi^q_{\lambda\lambda'}$, in the
sense that 
$\rho_{\lambda'} = \rho_\lambda
\circ
\Phi^q_{\lambda\lambda'}$. This
does not make sense as stated
because
$\Phi^q_{\lambda\lambda'}$ is
valued in the fraction algebra 
$\widehat{\mathcal T}^q_\lambda$
and not just in $\mathcal
T^q_\lambda$. A natural  approach
would be, for each $X\in 
\mathcal T^q_{\lambda'}$, to
write  the
rational
fraction $\Phi^q_{\lambda\lambda'}
(X)$ as the quotient
$PQ^{-1}$ of two polynomials $P$,
$Q\in
\mathcal T^q_\lambda $ and to
require that $\rho_{\lambda'} (X)
= \rho_\lambda(P)
\rho_\lambda(Q)^{-1}$. This of
course requires
$\rho(Q)$ to be invertible in
$\mathrm{End}(V)$, which creates
many problems in making the
definition consistent.  Actually,
for a general isomorphism 
$\Phi\col
\widehat{\mathcal T}^q_{\lambda'}
\rightarrow \widehat{\mathcal
T}^q_\lambda$ and for a
representation $\rho_\lambda\col 
\mathcal T^q_\lambda
\rightarrow
\mathrm{End} (V)$, it is
surprisingly difficult to determine
under which conditions on $\Phi$
and
$\rho_\lambda$ they define a
representation $\rho_\lambda
\circ \Phi \col \mathcal
T^q_{\lambda'}
\rightarrow
\mathrm{End} (V)$ in the above
sense. A lot of these problems
can be traced back to the fact
that, when adding up fractions
$PQ^{-1}$, the usual
technique of reduction to a common
denominator is much more
complicated in the
non-commutative context. 

We will use an
\emph{ad hoc} definition which
strongly uses the definition of
$\Phi_{\lambda\lambda'}^q$. 
After much work and provided we
consider all ideal triangulations
at the same time, it will
eventually turn out to be
equivalent to the above
definition. 

Given two ideal triangulations
$\lambda$ and $\lambda'$ and two
finite-dimensional
representations 
$\rho_\lambda\col \mathcal
T^q_\lambda
\rightarrow
\mathrm{End} (V)$ and 
$\rho_{\lambda'}\col \mathcal
T^q_{\lambda'}
\rightarrow
\mathrm{End} (V)$ of the
associated Chekhov-Fock algebras,
we say that $\rho_{\lambda'}$ is
\emph{compatible with}
$\rho_\lambda$ and we write 
 $\rho_{\lambda'} = \rho_\lambda
\circ
\Phi^q_{\lambda\lambda'}$ if, for
every generator
$X_i\in \mathcal T^q_{\lambda'}$,
we can write the rational fraction
$\Phi^q_{\lambda\lambda'}(X_i)
\in 
\widehat{\mathcal T}^q_\lambda$
as the quotient $P_iQ_i^{-1}$ of two
polynomials
$P_i$,
$Q_i\in\mathcal T^q_\lambda$ in
such a way that
$\rho_\lambda(Q_i)$ is invertible
in $\mathrm{End}(V)$ and
$\rho_{\lambda'}(X_i)=\rho_\lambda
(P_i)\rho_\lambda(Q_i)^{-1}$.
Note that $\rho_\lambda(P_i)$
then is also invertible by
consideration of
$\rho_{\lambda'}(X_i^{-1})$

At this point, it is not even
clear that the relation ``is
compatible with'' is symmetric
and transitive. A version of
these properties is provided by
the following lemma.

\begin{lem}
\label{lem:CriterionRepPolCore}
Consider a sequence of ideal
triangulations $\lambda_1$,
$\lambda_2$, \dots,
$\lambda_m$ and finite-dimensional
representations
$\rho_{\lambda_k}\col
\mathcal T^q_{\lambda_k}
\rightarrow
\mathrm{End} (V)$ such that  each
$\lambda_{k+1}$ is obtained from
$\lambda_k$ by a re-indexing or a
diagonal exchange. If in addition
$\rho_{\lambda_k} =
\rho_{\lambda_{k+1}} \circ
\Phi^q_{\lambda_{k+1}\lambda_k}$
for every $k$, then 
$\rho_{\lambda_1} =
\rho_{\lambda_{m}} \circ
\Phi^q_{\lambda_{m}\lambda_1}$
and $\rho_{\lambda_m} =
\rho_{\lambda_{1}} \circ
\Phi^q_{\lambda_{1}\lambda_m}$. 
\end{lem}

\begin{proof} We will prove that
$\rho_{\lambda_1}=\rho_{\lambda_m}
\circ \Phi_{\lambda_m\lambda_1}$
by induction on $m$. For this
purpose, assume the property true
for $m-1$. We need to show that,
for every generator $X_i$ of
$\mathcal T_{\lambda_1}^q$,
$\Phi_{\lambda_m\lambda_1}(X_i)$
can be written as a quotient
$PQ^{-1}$ where $P$, $Q\in
\mathcal T_{\lambda_1}^q$ are
such that $\rho_{\lambda_m}(P)$
and
$\rho_{\lambda_m}(Q)$ are
invertible and
$\rho_{\lambda_1}(X_i) =
\rho_{\lambda_m}(P)
\rho_{\lambda_m}(Q)^{-1}$. 

If $\lambda_m$ is obtained from
$\lambda_{m-1}$ by re-indexing,
then the property immediately
follows from the induction
hypothesis after re-indexing of
the $X_i$.

We can therefore restrict
attention to the case where
$\lambda_m$ is obtained from
$\lambda_{m-1}$ by one diagonal
exchange, along the
$i_0$--th component of
$\lambda_{m-1}$, say.

The general strategy of the proof
is fairly straightforward, but
the non-commutative context makes
it hard to control which elements
have an invertible image under
$\rho_{\lambda_m}$; this
requires more care that one might
have anticipated at first glance.

We need to be a little careful in
our notation. Let $\mathbb C\{
Z_1^{\pm1}, Z_2^{\pm1},
\dots, Z_n^{\pm1}\}$ denote the
algebra of non-commutative 
polynomials in the $2n$ variables
$Z_1$, $Z_2$, \dots, $Z_n$, 
$Z_1^{-1}$, $Z_2^{-1}$, \dots,
$Z_n^{-1}$. Given such a
polynomial
$P\in
\mathbb C\{ Z_1^{\pm1},
Z_2^{\pm1},
\dots, Z_n^{\pm1}\}$ and
invertible elements $A_1$, $A_2$,
\dots $A_n$ of an algebra
$\mathcal A$, we will denote by
$P(A_1, A_2, \dots, A_n)$ the
element of $\mathcal A$ defined
by replacing each $Z_i$ by the
corresponding $A_i$ and each
$Z_i^{-1}$ by $A_i^{-1}$. 

Consider the generator
$X_i \in \mathcal
T_{\lambda_1}^q$. By induction
hypothesis,
\begin{equation*}
\Phi^q_{\lambda_{m-1}\lambda_1}
(X_i) = P(X_1,  \dots, X_n)\,\,
Q(X_1,  \dots, X_n)^{-1}
\end{equation*} in
$\widehat{\mathcal
T}_{\lambda_{m-1}}^q$, for some 
non-commutative  polynomials
$P$ and $Q$ with \break
$\rho_{\lambda_{m-1}}(P(X_1, X_2,
\dots, X_n))$ and
$\rho_{\lambda_{m-1}}(Q(X_1,
X_2,
\dots, X_n))$ invertible in $
\mathrm{End}(V)$; beware that
$X_i$ represents a generator of
$\mathcal T_{\lambda_1}^q$ in the
left hand side of the equation,
and a generator of $\mathcal
T_{\lambda_{m-1}}^q$ in the right
hand side. In addition, 
\begin{equation*}
\rho_{\lambda_1} (X_i) =
\rho_{\lambda_{m-1}}(P(X_1, X_2,
\dots, X_n))\,\,
\rho_{\lambda_{m-1}}(Q(X_1, X_2,
\dots, X_n))^{-1}.
\end{equation*}

Then,
\begin{equation*}
\begin{split}
\Phi^q_{\lambda_m\lambda_1} (X_i)
&=
\Phi^q_{\lambda_m\lambda_{m-1}}
\circ
\Phi^q_{\lambda_{m-1}\lambda_1}
(X_i)\\ &=
\Phi^q_{\lambda_m\lambda_{m-1}}
(P(X_1,  \dots, X_n)) \,\,
\Phi^q_{\lambda_m\lambda_{m-1}}
(Q(X_1,  \dots, X_n))^{-1}\\ &=
P(\Phi^q_{\lambda_m\lambda_{m-1}}
(X_1), 
\dots,
\Phi^q_{\lambda_m\lambda_{m-1}}
(X_n))\\ &\quad\quad\quad
\quad\quad\quad\quad\quad
Q(\Phi^q_{\lambda_m\lambda_{m-1}}
(X_1), 
\dots,
\Phi^q_{\lambda_m\lambda_{m-1}}
(X_n))^{-1}.
\end{split}
\end{equation*} We are now facing
the problem of reducing these
quantities to a common
denominator, while controlling
the invertibility of the images
of denominators under
$\rho_{\lambda_m}$.

The ideal triangulation
$\lambda_m$ is obtained from
$\lambda_{m-1}$ by a diagonal
exchange along its
$i_0$--th component. By
inspection in the formulas
defining
$\Phi_{\lambda_m\lambda_{m-1}}$,
it follows that
$P(\Phi^q_{\lambda_m\lambda_{m-1}}
(X_1), 
\dots,
\Phi^q_{\lambda_m\lambda_{m-1}}
(X_n))$ is a polynomial in the
terms $X_j^{\pm 1}$, $(1+q
X_{i_0}^{\pm 1})^{-1}$ and
possibly $(1+q^3 X_{i_0}^{\pm
1})^{-1}$. In addition, whenever
a factor  $(1+q X_{i_0}^{\pm
1})^{-1}$ or $(1+q^3 X_{i_0}^{\pm
1})^{-1}$ appears, it is through
a relation such as
\begin{equation*}
\begin{split}
\Phi_{\lambda_m\lambda_{m-1}}
(X_j^{-1}) &=
X_j^{-1}(1+qX_{i_0})^{-1}\\
\text{or }\qquad
\Phi_{\lambda_m\lambda_{m-1}}
(X_j) &= (1+qX_{i_0}^{-1})^{-1}
(1+q^3X_{i_0}^{-1})^{-1}X_j
\end{split}
\end{equation*} (there are two
more possibilities), which
respectively give 
\begin{equation*}
\begin{split} &\rho_{\lambda_m}
(1+qX_{i_0}) =
\rho_{\lambda_{m-1}}(X_j)
\rho_{\lambda_m}(X_j^{-1}),\\
&\rho_{\lambda_m}
(1+q^3X_{i_0}^{-1})
\rho_{\lambda_m}
(1+qX_{i_0}^{-1}) =
\rho_{\lambda_m}(X_j)
\rho_{\lambda_{m-1}}(X_j^{-1}),
\end{split}
\end{equation*} or two more
relations, using the property that
$\rho_{\lambda_{m-1}} =
\rho_{\lambda_m} \circ
\Phi^q_{\lambda_m\lambda_{m-1}}$.
Since
$\rho_{\lambda_m} (X_j^{\pm1})$
and
$\rho_{\lambda_{m-1}}(X_j^{\pm1})$
are invertible and since $V$ is
finite-dimensional we conclude
that, for every
$(1+q X_{i_0}^{\pm 1})^{-1}$ or
$(1+q^3 X_{i_0}^{\pm 1})^{-1}$
appearing in
$P(\Phi^q_{\lambda_m\lambda_{m-1}}
(X_1), 
\dots,
\Phi^q_{\lambda_m\lambda_{m-1}}
(X_n))$, the corresponding element
$\rho_{\lambda_m}(1+q X_{i_0}^{\pm
1})$ or $\rho_{\lambda_m}(1+q^3
X_{i_0}^{\pm 1})$ is invertible
in $\mathrm{End}(V)$. 

Now, using the skew-commutativity
relations
\begin{equation*}
(1+q^{2k+1}X_{i_0}^{\pm1})X_j =X_j
(1+q^{2k
\pm
\sigma_{i_0 j}+1}X_{i_0}^{\pm1}) ,
\end{equation*} we can push all
the 
$(1+q^{2k+1}X_{i_0}^{\pm1})^{-1}$
to the right in the expression of \break
$P(\Phi^q_{\lambda_m\lambda_{m-1}}
(X_1), 
\dots,
\Phi^q_{\lambda_m\lambda_{m-1}}
(X_n))$, leading to a relation
\begin{equation*}
P(\Phi^q_{\lambda_m\lambda_{m-1}}
(X_1), 
\dots,
\Phi^q_{\lambda_m \lambda_{m-1}}
(X_n)) = P'(X_1, \dots, X_n)
R(X_{i_0})^{-1}
\end{equation*}  where $ P'(X_1,
\dots, X_n)$ is a Laurent
polynomial in the
$X_j$ and where
$R(X_{i_0})$ is a 1--variable
Laurent polynomial product of
terms
$(1+q^{2k+1}X_{i_0}^{\pm1})$. In
addition, applying
$\rho_{\lambda_m}$ to both sides
of the above skew-commutativity
relation, we see that
$\rho_{\lambda_m}
(1+q^{2k+1}X_{i_0}^{\pm1})$ is
invertible in $\mathrm{End}(V)$
whenever a term
$(1+q^{2k+1}X_{i_0}^{\pm1})^{-1}$
appears in this process.
Therefore,
$\rho_{\lambda_m}(R(X_{i_0}))$ is
invertible. 

We will now perform essentially
the same computations in
$\mathrm{End}(V)$. Since
$\rho_{\lambda_{m-1}} =
\rho_{\lambda_m} \circ
\Phi^q_{\lambda_m\lambda_{m-1}}$,
\begin{equation*}
\begin{split}
\rho_{\lambda_{m-1}}(P(X_1, 
\dots, X_n)) &= 
P(\rho_{\lambda_{m-1}}(X_1), 
\dots,
\rho_{\lambda_{m-1}}(X_n))\\ &= 
P(\rho_{\lambda_m}
\circ\Phi^q_{\lambda_m
\lambda_{m-1}}(X_1), 
\dots, \rho_{\lambda_m} \circ
\Phi^q_{\lambda_m\lambda_{m-1}}(X_n))
\end{split}
\end{equation*} The same
manipulations as above, but
replacing the $X_j$ by the
$\rho_{\lambda_m} (X_j) \in
\mathrm{End}(V)$ (which satisfy
the same relations),  yield
\begin{equation*}
\begin{split}
\rho_{\lambda_{m-1}}(P(X_1, 
\dots, X_n))  &= 
P(\rho_{\lambda_m}
\circ\Phi^q_{\lambda_m
\lambda_{m-1}} (X_1), 
\dots, \rho_{\lambda_m} \circ
\Phi^q_{\lambda_m
\lambda_{m-1}}(X_n))
\\ &= P'( \rho_{\lambda_m}(X_1),
\dots,
\rho_{\lambda_m} (X_n) )\,
R(\rho_{\lambda_m}
(X_{i_0}))^{-1}\\ &=
\rho_{\lambda_m} (P'( X_1, \dots,
X_n ))\,
\rho_{\lambda_m} (R(X_{i_0}))^{-1}
\end{split}
\end{equation*} In particular,
since 
$\rho_{\lambda_{m-1}}(P(X_1, X_2,
\dots, X_n))$ is invertible by
definition of $P$ and
$Q$ and since
$\rho_{\lambda_m} (R(X_{i_0}))$ is
invertible by construction, we
conclude that
$\rho_{\lambda_m} (P'( X_1,
\dots, X_n ))$ is invertible. 

Similarly, we can write 
\begin{equation*}
Q(\Phi^q_{\lambda_m\lambda_{m-1}}
(X_1), 
\dots,
\Phi^q_{\lambda_m\lambda_{m-1}}
(X_n)) = Q'(X_1, \dots, X_n)
S(X_{i_0})^{-1}
\end{equation*}  for some Laurent
polynomials $Q'(X_1, \dots, X_n)$
and
$S(X_{i_0})$, in such a way that 
$\rho_{\lambda_m} (Q_i( X_1,
\dots, X_n ))$ and
$\rho_{\lambda_m} (S(X_{i_0}))$
are invertible in
$\mathrm{End}(V)$, and
\begin{equation*}
\rho_{\lambda_{m-1}}(Q(X_1, 
\dots, X_n))  = \rho_{\lambda_m}
(Q'( X_1,
\dots, X_n ))\,
\rho_{\lambda_m}(S(X_{i_0}))^{-1}.
\end{equation*}

We are now ready to conclude.
Indeed, we showed that
\begin{equation*}
\begin{split}
\Phi^q_{\lambda_m\lambda_1} (X_i) 
&=
P(\Phi^q_{\lambda_m\lambda_{m-1}}
(X_1), 
\dots,
\Phi^q_{\lambda_m\lambda_{m-1}}
(X_n))\\
&\quad\quad\quad\quad\quad\quad
\quad\quad\quad
Q(\Phi^q_{\lambda_m\lambda_{m-1}}
(X_1), 
\dots,
\Phi^q_{\lambda_m\lambda_{m-1}}
(X_n))^{-1}\\ &=
\left( P'(X_1, \dots, X_n)
R(X_{i_0})^{-1}
\right)
\left( Q'(X_1, \dots, X_n)
S(X_{i_0})^{-1}
\right)^{-1}\\ &=
\left( P'(X_1, \dots, X_n)
S(X_{i_0})
\right)
\left( Q'(X_1, \dots, X_n)
R(X_{i_0})
\right)^{-1}
\end{split}
\end{equation*} since
$R(X_{i_0})$ and
$S(X_{i_0})$ commute. Similarly, 
\begin{equation*}
\begin{split}
\rho_{\lambda_1} (X_i) &=
\rho_{\lambda_{m-1}}(P(X_1, X_2,
\dots, X_n))
\rho_{\lambda_{m-1}}(Q(X_1, X_2,
\dots, X_n))^{-1} \\ &=
\left(
\rho_{\lambda_m} (P'( X_1, \dots,
X_n ))\,
\rho_{\lambda_m} (R(X_{i_0}))^{-1}
\right)\\
&\quad\quad\quad\quad\quad\quad
\quad\quad
\left(
\rho_{\lambda_m}(Q'( X_1, \dots,
X_n ))\,
\rho_{\lambda_m}(S(X_{i_0}))^{-1}
\right)^{-1}\\ &=
\left(
\rho_{\lambda_m} (P'( X_1, \dots,
X_n ))\,
\rho_{\lambda_m}(S(X_{i_0}))
\right)\\
&\quad\quad\quad\quad\quad\quad
\left(
\rho_{\lambda_m}(Q'( X_1, \dots,
X_n ))\,
\rho_{\lambda_m}(R(X_{i_0}))
\right)^{-1}\\ &=
\rho_{\lambda_m}(P'( X_1, \dots,
X_n )S(X_{i_0}))\,\,
\rho_{\lambda_m} (Q'( X_1, \dots,
X_n )R(X_{i_0}))^{-1}
\end{split}
\end{equation*} By definition,
this means that
$\rho_{\lambda_1}=\rho_{\lambda_m}
\circ \Phi_{\lambda_m\lambda_1}$,
as desired. 

There remains to prove the second
statement that $\rho_{\lambda_m} =
\rho_{\lambda_{1}} \circ
\Phi^q_{\lambda_{1}\lambda_m}$.
For this, note that the property
that
$\rho_{\lambda_k} =
\rho_{\lambda_{k+1}} \circ
\Phi^q_{\lambda_{k+1}\lambda_k}$
implies that $\rho_{\lambda_{k+1}}
=
\rho_{\lambda_{k}} \circ
\Phi^q_{\lambda_{k}\lambda_{k+1}}$
for every $k$, using the explicit
form of
$\Phi^q_{\lambda_{k+1}\lambda_k}$
and
$\Phi^q
_{\lambda_{k}\lambda_{k+1}}$ as
well as  arguments which are
similar to and much simpler than
the ones we just used.   The
property that
$\rho_{\lambda_m} =
\rho_{\lambda_{1}} \circ
\Phi^q_{\lambda_{1}\lambda_m}$
then immediately follows by
symmetry.
\end{proof}

A \emph{representation of the
polynomial core} ${\mathcal
T}^q_S$ over the
 vector space
$V$ is a family of
representations 
$\rho_\lambda\col \mathcal
T^q_\lambda
\rightarrow
\mathrm{End} (V)$ defined for
each ideal triangulation
$\lambda \in \Lambda(S)$, such
that any two $\rho_{\lambda'}$ and
$\rho_\lambda$ are compatible in
the above sense.  Lemma~\ref
{lem:CriterionRepPolCore} shows
that it suffices to check this
condition on pairs of ideal
triangulations which are obtained
from each other by one re-indexing
or one diagonal exchange.  We will
see in the next sections that the
polynomial core admits many
representations.

Before closing this section, we
indicate the following result,
which shows that our
definition of compatibility
coincides with the condition we
had in mind at the beginning of
this section. 

\begin{lem} 
\label{lem:RepPolCore}
Let $\rho = 
\{\rho_\lambda\col
\mathcal T^q_\lambda
\rightarrow
\mathrm{End} (V)\}_{\lambda \in
\Lambda(S)}$ be a
finite-dimensional irreducible
representation of the polynomial
core $\mathcal T^q_S$ of the
quantum Teichm\"uller space
$\widehat{\mathcal T}^q_S$. Then,
for every $X' \in \mathcal
T_{\lambda'}^q$, its image
$\Phi^q_{\lambda\lambda'}(X')
\in \widehat{\mathcal
T}_{\lambda}^q$ can be written as
$\Phi^q_{\lambda\lambda'}(X')
=PQ^{-1} =(Q')^{-1}P'$ with
$P$, $Q\in \mathcal T_\lambda^q$
and with $\rho_\lambda(Q)$ and
$\rho_\lambda(Q')$ invertible in
$\mathrm {End}(V)$. In addition,
for any such decomposition of 
$\Phi^q_{\lambda\lambda'}(X')$, 
$\rho_{\lambda'}(X')$ is then
equal to
$\rho_\lambda(P)\rho(Q)^{-1} =
\rho(Q')^{-1} \rho_\lambda(P')$.
\qed
\end{lem}
\begin{proof}
This is proved by
arguments almost identical to the
ones we used for
Lemma~\ref
{lem:CriterionRepPolCore}, by
induction on the number of
diagonal exchanges needed to go
from $\lambda$ to $\lambda'$.
However, it is worth mentioning
that the easy algebraic
manipulation leading to the last
statement simultaneously uses the
left and right decompositions
$PQ^{-1}$ and $(Q')^{-1}P'$ of 
$\Phi^q_{\lambda\lambda'}(X')$. 
\end{proof}

\section{The non-quantum shadow
of a representation}
\label{sect:NonQuantumShadow}

By
Theorems~\ref{thm:ClassRepCFodd}
and \ref{thm:ClassRepCFeven}, an
irreducible finite-dimensional 
representation 
$\rho_\lambda\col \mathcal
T^q_\lambda
\rightarrow
\mathrm{End} (V)$ of the
Chekhov-Fock algebra is
classified, up to a finite number
of choices of certain roots,  by
numbers
$x_i
\in
\mathbb C^*$ associated to the
components
$\lambda_i$ of $\lambda$. By 
Theorem~\ref{thm:ClassRepCFodd}
or by inspection, the same
numbers $x_i$ completely
determine a representation
$\rho_\lambda^1\col \mathcal
T^1_\lambda
\rightarrow
\mathrm{End} (\mathbb C)$ of the
commutative algebra $\mathcal
T^1_\lambda$ corresponding to the
non-quantum (also called
semi-classical in the physics
literature) case where
$q=1$. We will say that
$\rho_\lambda^1$ is the
\emph{non-quantum shadow}, or the
\emph{semi-classical shadow}, of
the representation
$\rho_\lambda$.

Interpreting the numbers $x_i
\in
\mathbb C^*$ as a non-quantum
representation
$\rho_\lambda^1\col \mathcal
T^1_\lambda
\rightarrow
\mathrm{End} (\mathbb C)$ may
sound really pedantic at first.
However, the remainder of this
paper hinges on the following
computation which shows that, for
a suitable choice of
$q$, the map
$\rho_\lambda
\mapsto
\rho_\lambda^1$  is well-behaved
with respect to the coordinate
changes
$\Phi^q_{\lambda\lambda'}$ and
$\Phi^1_{\lambda\lambda'}$.

\begin{lem}
\label{lem:QuantumNonQuantum} Let
$q$ be such that
$q^2$ is a primitive $N$--th root
of unity and such that
$q^N=(-1)^{N+1}$ (for instance 
$q=-
\mathrm e^{\pi \mathrm i/N}$).
Suppose that the two ideal
triangulations
$\lambda$ and $\lambda'$ of the
surface $S$ are obtained from
each other by a diagonal exchange
or by a re-indexing, and consider
two irreducible finite-dimensional
representations 
$\rho_\lambda\col \mathcal
T^q_\lambda
\rightarrow
\mathrm{End} (V)$ and 
$\rho_{\lambda'}\col \mathcal
T^q_{\lambda'}
\rightarrow
\mathrm{End} (V)$ such that 
$\rho_{\lambda'} = \rho_{\lambda}
\circ
\Phi^q_{\lambda\lambda'}$ in the
sense of \S\ref{sect:PolCore}. If 
$\rho_{\lambda}^1\col \mathcal
T^1_\lambda
\rightarrow
\mathrm{End} (\mathbb C)$ and 
$\rho_{\lambda'}^1\col \mathcal
T^1_{\lambda'}
\rightarrow
\mathrm{End} (\mathbb C)$ are the
respective non-quantum shadows of
$\rho_{\lambda}$ and
$\rho_{\lambda'}$, then
$\rho_{\lambda'}^1 =
\rho_{\lambda}^1
\circ
\Phi^1_{\lambda\lambda'}$.
\end{lem}

\begin{proof} Recall  that
$\rho_{\lambda}^1$ is determined
by the property that
$\rho_{\lambda}^1(X_i) = x_i \in
\mathbb C^* \subset
\mathrm{End}(\mathbb C)$, where
$x_i$ is the number such that
$\rho_{\lambda}(X_i^N) = x_i\,
\mathrm{Id}_V$. Similarly, 
$\rho_{\lambda'}^1(X_i) = x_i'$
where
$x_i'$ is such that
$\rho_{\lambda'}(X_i^N)
= \rho_{\lambda}\circ
\Phi_{\lambda\lambda'}
(X_i^N)=x_i'\,
\mathrm{Id}_V$. In particular,
the property is immediate when
$\lambda'$ is obtained from
$\lambda$ by a re-indexing of its
components. 

Suppose that $\lambda'$ is
obtained from $\lambda$ by an
embedded $i$--th diagonal
exchange. Label the four sides of
the square $Q$ supporting the
exchange counterclockwise as
$\lambda_j$, $\lambda_k$,
$\lambda_l$ and $\lambda_m$, in
such a way that the diagonal
$\lambda_i$ goes from the
$\lambda_j\lambda_k$ corner to
the $\lambda_l\lambda_m$ corner,
as in Figure~\ref{pict:DiagEx}. 

By definition of
$\Phi_{\lambda\lambda'}^q$,
$\Phi_{\lambda\lambda'}^q(X_i^N) =
X_i^{-N}$. Using Lemma~\ref
{lem:RepPolCore}, it
follows that
$\rho_{\lambda'}(X_i^N) =
\rho_{\lambda}(X_i^N)^{-1}$, so
that
$x_i'=x_i^{-1}$. 

Because $X_jX_i=q^2X_iX_j$,  the
quantum binomial formula (see for
instance
\cite[\S IV.2]{Kas}) shows
that 
\begin{equation*}
\begin{split}
\Phi_{\lambda\lambda'}^q(X_j^N) 
&=
\Phi_{\lambda\lambda'}^q(X_j)^N
 =(X_j + q X_iX_j)^N \\ &= X_j^N
+ (qX_iX_j)^N = X_j^N + q^N
q^{N(N-1)} X_i^N X_j^N \\ &=
X_j^N + X_i^N X_j^N.
\end{split}
\end{equation*} Indeed, most of
the quantum binomial coefficients
are 0 since $q^2$ is a primitive
$N$--th root of unity. Note that
we also used our hypothesis that
$q^N=(-1)^{N+1}$ for the last
equality.  It follows that $x'_j=
x_j + x_ix_j = (1+x_i)x_j$.

To compute $x_k$, it is easier to
consider
\begin{equation*}
\begin{split}
\Phi_{\lambda\lambda'}^q
(X_k^{-N})  &=
\Phi_{\lambda\lambda'}^q
(X_k)^{-N}
 =(X_k^{-1} + q
X_k^{-1}X_i^{-1})^N \\ &= X_k^{-N}
+ (qX_k^{-1}X_i^{-1})^N =
X_k^{-N} + q^N q^{N(N-1)} X_k^{-N}
X_i^{-N} \\ &= X_k^{-N} + X_k^{-N}
X_i^{-N}.
\end{split}
\end{equation*} 
Applying Lemma~\ref
{lem:RepPolCore}, we conclude that
$x'_k = (x_k^{-1} +
x_k^{-1}x_i^{-1})^{-1}
=(1+x_i^{-1})^{-1}x_k$.

Similar computations hold for
$x_l'$ and $x_m'$. We conclude
that  $x_i'=x_i^{-1}$,
$x'_j=  (1+x_i)x_j$, $x'_k
=(1+x_i^{-1})^{-1}x_k$,  $x'_l=
(1+x_i)x_l$, $x'_m 
=(1+x_i^{-1})^{-1}x_m$ and
$x_h'=x_h$ if $h\not\in
\{i,j,k,l,m\}$. By definition of
$\Phi^1_{\lambda\lambda'}$, this
just means that
$\rho_{\lambda'}^1 =
\rho_{\lambda}^1
\circ
\Phi^1_{\lambda\lambda'}$.

This completes the proof for an
embedded diagonal exchange.

We now consider  non-embedded
diagonal exchanges. Keeping the
same labelling conventions as
before, suppose that we are in
the case called Case~2 earlier,
namely where
$\lambda_j=\lambda_l$ and
$\lambda_k \not= \lambda_m$. In
this situation, $X_jX_i=q^4
X_iX_j$ in $\mathcal
T^q_\lambda$, which obliges us to
use different arguments according
to the parity of $N$.

If $N$ is odd, then $q^4$ is
still a primitive $N$--th root of
unity, and the quantum binomial
formula again shows that
\begin{equation*}
\begin{split}
\Phi_{\lambda\lambda'}^q(X_j^N)  &
=\bigl( (1 + q X_i) (1 + q^3
X_i)X_j
\bigr)^N =\bigl( U + q X_i U
\bigr)^N
\\ &= U^N + (qX_iU)^N = U^N + q^N
q^{2N(N-1)} X_i^N U^N \\ &= (1 +
X_i^N)U^N
\end{split}
\end{equation*} where $U=(1 + q^3
X_i)X_j$; note for this that
$UX_i=q^4X_iU$, and also use
$q^N=(-1)^{N+1}=1$. Another
application of the quantum
binomial
 formula gives
\begin{equation*}
\begin{split} U^N  & =\bigl( X_j
+ q^3 X_iX_j
\bigr)^N =X_j^N + (q^3X_iX_j)^N
\\ &=  X_j^N + q^{3N} q^{2N(N-1)}
X_i^N X_j^N = (1 +  X_i^N)X_j^N
\end{split}
\end{equation*} so that
$\Phi_{\lambda\lambda'}(X_j^N)= (1
+  X_i^N)^2 X_j^N$. This implies
that $x'_j=(1+x_i)^2 x_j$. The
same computations as in the
embedded diagonal exchange case
give $x_i'=x_i^{-1}$, 
$x'_k =(1+x_i^{-1})^{-1}x_k$,  
$x'_m  =(1+x_i^{-1})^{-1}x_m$ and
$x_h'=x_h$ if $h\not\in
\{i,j,k,l,m\}$. By definition of
$\Phi^1_{\lambda\lambda'}$, this
implies that
$\rho_1' = \rho_1
\circ
\Phi^1_{\lambda\lambda'}$ in this
case as well.

When $N$ is even, there is a new
twist because $q^4$ is now a
primitive
$\frac N2$--th root of unity. For
$U$ as above, the quantum binomial
formula gives in this case
\begin{equation*}
\begin{split}
\Phi_{\lambda\lambda'}^q\bigl
(X_j^{\frac N2} \bigr)  & =\bigl(
(1 + q X_i) (1 + q^3 X_i)X_j
\bigr)^{\frac N2} =\bigl( U + q
X_i U
\bigr)^{\frac N2}
\\ &= U^{\frac N2} +
(qX_iU)^{\frac N2} = U^{\frac N2}
+ q^{\frac N2} q^{\frac
{N(N-2)}2} X_i^{\frac N2} U^{\frac
N2} \\ &= \bigl(1 + (-1)^{\frac
{N-2}2 } q^{\frac N2} X_i^{\frac
N2} \bigr) U^{\frac N2}
\end{split}
\end{equation*} and
\begin{equation*}
\begin{split} U^{\frac N2}  &
=\bigl( X_j + q^3 X_iX_j
\bigr)^{\frac N2} =X_j^{\frac N2}
+ (q^3X_iX_j)^{\frac N2}
\\ &=  X_j^{\frac N2} +
q^{\frac{3N}2} q^{\frac
{N(N-2)}2} X_i^{\frac N2}
X_j^{\frac N2} = \bigl(1 +
(-1)^{\frac {N}2 } q^{\frac N2}
X_i^{\frac N2} \bigr) X_j^{\frac
N2},
\end{split}
\end{equation*} using the fact
that $q^N = (-1)^{N+1}=-1$. It
follows that 
$
\Phi_{\lambda\lambda'}^q\bigl
(X_j^{\frac N2} \bigr) = \bigl(
1- q^N X_i^N 
\bigr) X_j^{\frac N2} = \bigl( 1+
X_i^N 
\bigr) X_j^{\frac N2} $. Noting
that $X_i^N$ and
$X_j^{\frac N2}$ commute, we
conclude that 
$\Phi_{\lambda\lambda'}^q
\bigl( X_j^N \bigr)=
\Phi_{\lambda\lambda'}^q
\bigl( X_j^{\frac N2} \bigr)^2= (1
+  X_i^N)^2 X_j^N$ in this case as
well. Therefore,
$x'_j=(1+x_i)^2 x_j$,
$x_i'=x_i^{-1}$, 
$x'_k =(1+x_i^{-1})^{-1}x_k$,  
$x'_m  =(1+x_i^{-1})^{-1}x_m$ and
$x_h'=x_h$ if $h\not\in
\{i,j,k,l,m\}$ as before. This
again implies that
$\rho_{\lambda'}^1 =
\rho_{\lambda'}^1
\circ
\Phi^1_{\lambda\lambda'}$ in this
case.

The remaining types of
non-embedded diagonal exchanges
are treated in the same way, using
the above computations.
\end{proof}

Note that the conditions that
$q^2$ is a primitive $N$--th root
of unity and
$q^N=(-1)^{N+1}$ are equivalent to
the property that $q$ is a
primitive $N$--th root of
$(-1)^{N+1}$, which is shorter to
state. The combination of
Lemmas~\ref{lem:QuantumNonQuantum}
and \ref{lem:CriterionRepPolCore}
immediately gives:

\begin{thm} 
\label{thm:NonQuantumShadow} Let
$q$ be a primitive $N$--th root
of $(-1)^{N+1}$. If  $\rho =
\{\rho_\lambda\col
\mathcal T^q_\lambda
\rightarrow
\mathrm{End} (V)\}_{\lambda \in
\Lambda(S)}$ is a
finite-dimensional irreducible
representation of the polynomial
core $\mathcal T^q_S$ of the
quantum Teichm\"uller space
$\widehat{\mathcal T}^q_S$, then
the non-quantum shadows of the
$\rho_\lambda$ form a
representation $\rho^1 =
\{\rho_\lambda^1\col
\mathcal T^1_\lambda
\rightarrow
\mathrm{End} (\mathbb
C)\}_{\lambda
\in
\Lambda(S)}$ of the non-quantum
polynomial core $\mathcal T^1_S$.
\qed
\end{thm}

We will say that the
representation $\rho^1$ of the
polynomial core $\mathcal T_S^1$
is the \emph{non-quantum shadow}
of the representation
$\rho$ of the polynomial 
core $\mathcal T^q_\lambda$. 

We now show that every
representation of the
non-quantum polynomial core
$\mathcal T^1_S$ is the shadow
of several representations of the
quantum polynomial core
$\mathcal T_S^q$.  

\begin{lem}
\label{lem:NonQuantumGivesQuantum}
Let the ideal triangulation
$\lambda'$ be obtained from
$\lambda$ by a re-indexing or by a
diagonal exchange. Consider an
irreducible finite-dimensional
representation
$\rho_\lambda \col
\mathcal T_\lambda^q \rightarrow
\mathrm{End}(V) $,
with non-quantum shadow
$\rho_\lambda^1 \col
\mathcal T_\lambda^1 \rightarrow
\mathrm{End}(\mathbb C)$. If
there exists a non-quantum
representation
$\rho_{\lambda'}^1\col
\mathcal T_{\lambda'}^1
\rightarrow
\mathrm{End}(\mathbb C)= \mathbb
C^*$ with
$\rho_{\lambda'}^1 =
\rho_\lambda^1
\circ \Phi_{\lambda\lambda'}^1$,
then there exists a unique
representation
$\rho_{\lambda'}\col
\mathcal T_{\lambda'}^q
\rightarrow
\mathrm{End}(V)$  with 
$\rho_{\lambda'} =
\rho_\lambda
\circ \Phi_{\lambda\lambda'}^q$
and with shadow
$\rho_{\lambda'}^1$.
\end{lem}

\begin{proof} The property is
immediate for a re-indexing.

Suppose that $\lambda'$ is
obtained from $\lambda$ by an
embedded diagonal exchange along
the component $\lambda_i$. Label
the components of $\lambda$
bounding the square $Q$ where the
diagonal exchange takes place as
$\lambda_j$,
$\lambda_k$, $\lambda_l$ and
$\lambda_m$, as in Figure~\ref
{pict:DiagEx}.
By inspection of the formulas
defining
$\Phi_{\lambda\lambda'}^q$,
$\rho_{\lambda'} (X_s^{\pm1})=
\rho_\lambda
\circ \Phi_{\lambda\lambda'}^q
(X_s^{\pm1})$ will be defined if
$\rho_\lambda (1 +qX_{i})$ and
$\rho_\lambda (1 +qX_{i}^{-1})$
are invertible in
$\mathrm{End}(V)$. As in the
proof of Lemma~\ref
{lem:QuantumNonQuantum}, 
\begin {equation*}
\begin{split}
\rho_\lambda \left(
(1+qX_i)X_j\right)^N &=
(1+\rho_\lambda(X_i^N))
\rho_\lambda(X_j^N)\\ &=
(1+\rho_\lambda^1(X_i))
\rho_\lambda^1 (X_j)\,
\mathrm{Id}_V\\ &= 
\rho_{\lambda'}^1 (X_j)\,
\mathrm{Id}_V.
\end{split}
\end{equation*} 
Since $\rho_{\lambda'}^1 (X_j)
\not=0$, it follows that
$\rho_\lambda \left(
(1+qX_i)X_j\right)$ is
invertible, and therefore so is
$\rho_\lambda \left(
(1+qX_i)\right)$. A similar
consideration of
$\rho_\lambda(X_k^{-1} (1
+qX_{i}^{-1}))^N$ proves the
invertibility of
$\rho_\lambda (1 +qX_{i}^{-1})$.

This defines $\rho_{\lambda'}$ on
the generators $X_s^{\pm1}$. By
inspection, it is compatible with
the skew-commutativity relations
$X_sX_t = q^{2\sigma_{st}}
X_tX_s$ and consequently extends
to an algebra homomorphism 
$\rho_{\lambda'}\col
\mathcal T_{\lambda'}^q
\rightarrow
\mathrm{End}(V)$. Its non-quantum
shadow is equal to
$\rho_{\lambda'}^1$. 

The case of a non-embedded
diagonal exchange is treated in
the same way, applying again the
computations of the proof of
Lemma~\ref
{lem:QuantumNonQuantum}. 
\end{proof}

\begin{thm}
\label{thm:NonQuantumGivesRep}
Let $q$ be a primitive $N$--th
root of $(-1)^{N+1}$. Up to
isomorphism, every representation
$\rho^1 =
\{\rho_\lambda^1\col
\mathcal T^1_\lambda
\rightarrow
\mathrm{End} (\mathbb
C)\}_{\lambda
\in
\Lambda(S)}$ of the non-quantum
polynomial core $\mathcal T_S^1$
is the non-quantum shadow of
exactly
$N^p$ if
$N$ is odd, and
$2^{2g}N^p$ is $N$ is even,
irreducible finite-dimensional 
representations
$\rho =
\{\rho_\lambda\col
\mathcal T^q_\lambda
\rightarrow
\mathrm{End} (V)\}_{\lambda \in
\Lambda(S)}$ of the polynomial
core $\mathcal T_\lambda^q$,
where $p$ is the number of
punctures of $S$ and $g$ is its
genus. 
\end{thm}

\begin{proof}  Fix an ideal
triangulation
$\lambda$. By
Theorem~\ref{thm:ClassRepCFodd}
or \ref{thm:ClassRepCFeven},
according to the parity of $N$,
there are $N^p$ or $2^{2g}N^p$
isomorphism classes of
irreducible finite-dimensional
representations
$\rho_\lambda\col
\mathcal T^q_\lambda
\rightarrow
\mathrm{End} (V)$ with
non-quantum shadow
$\rho_\lambda^1$. The combination
of Lemmas~\ref
{lem:CriterionRepPolCore} and
\ref {lem:NonQuantumGivesQuantum}
shows that each such
representation
$\rho_\lambda$ uniquely extends
to a representation of the
polynomial core $\mathcal T_S^q$. 
\end{proof}

\section{Pleated surfaces and the
hyperbolic shadow of a
representation}
\label{sect:HyperShadow}

We have just showed that the
representation theory of the
polynomial core $\mathcal T_S^q$
is, up to finitely many 
ambiguities, controlled by the
representation theory of the
non-quantum polynomial core
$\mathcal T_S^1$. It is now time
to remember that the non-quantum
coordinate changes
$\Phi_{\lambda\lambda'}^1$ were
specially designed to mimic the
coordinate changes between shear
coordinates  for the
Teichm\"uller space of the
surface $S$, or more precisely
for the enhanced Teichm\"uller
space as defined in \cite{Liu}.
We are going to take advantage of
this geometric context. 

However, when considering the
weights associated to a
non-quantum representation, we
subreptitiously moved from  real
to complex numbers. This leads us
to consider the complexification
of the Teichm\"uller space, when
considered as a real analytic
manifold. This complexification
has a nice geometric
interpretation, based on the fact
that the complexification of the
orientation-preserving isometry
group
$\mathrm{PSL}_2(\mathbb R)$ of
the hyperbolic plane
$\mathbb H^2$ is the
orientation-preserving isometry
group
$\mathrm{PSL}_2(\mathbb C)$ of
the hyperbolic 3--space
$\mathbb H^3$.  For this, we will
use the technical tool of pleated
surfaces, which is now classical
in 3--dimensional hyperbolic
geometry \cite{Th1, CEG, Bon}.

Let $\lambda$ be an ideal
triangulation of the surface $S$. 
A
\emph{pleated surface} with
\emph{pleating locus} $\lambda$ is
a pair
$\bigl( \widetilde f , r\bigr)$,
where $\widetilde f \col
\widetilde S \rightarrow \mathbb
H^3$ is a map from the universal
covering
$\widetilde S$ of $S$ to the
hyperbolic 3--space $\mathbb
H^3$, and where $r\col \pi_1(S)
\rightarrow
\mathrm{PSL}_2(\mathbb C)$ is a
group homomorphism from the
fundamental group of
$S$ to the group of
orientation-preserving isometries
of
$\mathbb H^3$, such that:
\begin{enumerate}
\item $\widetilde f$
homeomorphically sends each
component of the preimage
$\widetilde \lambda$ of $\lambda$
to a complete geodesic of $\mathbb
H^3$;
\item $\widetilde f$
homeomorphically sends the
closure of each component of
$\widetilde S-
\widetilde \lambda$ to an ideal
triangle in $\mathbb H^3$, namely
one whose three vertices are on
the sphere at infinity
$\partial_\infty \mathbb H^3$ of
$\mathbb H^3$;
\item $\widetilde f$ is
$r$--equivariant in the sense
that $\widetilde f(\gamma
\widetilde x) = r(\gamma)
\widetilde f(
\widetilde x)$ for every
$\widetilde x \in \widetilde S$
and $\gamma \in \pi_1(S)$. 
\end{enumerate}

In classical examples arising
from geometry, the homomorphism 
$r$ has discrete image, so that
$\widetilde f$ induces a map
$f\col S
\rightarrow \mathbb H^3/
r(\pi_1(S))$ to the quotient
orbifold $\mathbb H^3/
r(\pi_1(S))$. The map $f$ is
totally geodesic on $S-\lambda$,
and is bent along a geodesic
ridge at the components of
$\lambda$.

The geometry of the pleated
surface $\bigl( \widetilde f ,
r\bigr)$ is completely described
by complex numbers $x_i \in
\mathbb C^*$ associated to the
components $\lambda_i$ as
follows. Consider the upper
half-space model for $\mathbb
H^3$, bounded by the Riemann
sphere
$\widehat{\mathbb C}=
\mathbb C \cup
\{\infty\}$. Arbitrarily orient
$\lambda_i$ and lift it to an
oriented component
$\widetilde\lambda_i$ of
$\widetilde\lambda$. Let
$T_{\mathrm l}$ be the component
of $\widetilde S
-\widetilde\lambda$ that is on
the left of
$\widetilde\lambda_i$, and let 
$T_{\mathrm r}$ be the component
on the right, defined with
respect to  the orientations of
$\widetilde\lambda_i$ and
$\widetilde S$. Let $z_+$ and
$z_-\in \widehat{\mathbb C}$ be
the positive and negative end
points of the oriented geodesic
$\widetilde f\bigl(
\widetilde\lambda_i \bigr)$ of
$\mathbb H^3$, let $z_{\mathrm
l}$ be the vertex of the ideal
triangle $\widetilde f(T_{\mathrm
l})$ that is different from
$z_\pm$ and, likewise, let
$z_{\mathrm r}$ be the third
vertex of $T_{\mathrm r}$. Then
$x_i$ is defined as the
cross-ratio
\begin{equation*} x_i =
-\frac{(z_{\mathrm l}-z_+)
(z_{\mathrm r}-z_-)}{(z_{\mathrm
l}-z_-) (z_{\mathrm r}-z_+)}.
\end{equation*} Note that $x_i$
is different from
$0$ and $\infty$, because the
vertex sets $\{z_+, z_-,
z_{\mathrm l} \}$ and $\{z_+, z_-,
z_{\mathrm r} \}$ of the ideal
triangles $\widetilde
f(T_{\mathrm l})$ and $\widetilde
f(T_{\mathrm r})$ each consist of
three distinct points. Also,
reversing the orientation of
$\lambda_i$ exchanges
$z_+$ and $z_-$, but also
exchanges $z_{\mathrm l}$ and
$z_{\mathrm r}$ so that $x_i$ is
unchanged. Similarly, $x_i$ is
independent of the choice of the
lift $\widetilde \lambda_i$ by
invariance of cross-ratios under
hyperbolic isometries. 

By definition, $x_i \in \mathbb
C^*$ is the \emph{exponential
shear-bend parameter} of the
pleated surface $\bigl(
\widetilde f , r\bigr)$  along the
component $\lambda_i$ of
$\lambda$.  Geometrically, the
imaginary part of $\log x_i$
(defined modulo
$2\pi \mathrm i$) is the external
dihedral angle of the ridge
formed by $\widetilde f \bigl(
\widetilde S \bigr)$ near the
preimage of $\lambda_i$. The real
part of $\log x_i$ is the oriented
distance from $z_{\mathrm l}'$
to  $z_{\mathrm r}'$ in the
oriented geodesic  $\widetilde
f\bigl( \widetilde\lambda_i
\bigr)$, where 
$z_{\mathrm l}'$ and  $z_{\mathrm
r}'$ are the respective
orthogonal  projections of
$z_{\mathrm l}$ and  $z_{\mathrm
r}$ to 
$\widetilde f\bigl(
\widetilde\lambda_i
\bigr)$. See for instance
\cite{Bon}. 

Two pleated surfaces $\bigl(
\widetilde f , r\bigr)$ and
$\bigl( \widetilde f' , r'\bigr)$
are
\emph{isometric} if there is a
hyperbolic isometry $A\in
\mathrm{PSL}_2(\mathbb C)$ and a
lift $\widetilde\varphi \col
\widetilde S
\rightarrow \widetilde S$ of an
isotopy of $S$ such that
$\widetilde f' = A \circ
\widetilde f \circ
\widetilde\varphi$ and $r'(\gamma)
= A \,r(\gamma) A^{-1}$ for every
$\gamma \in \pi_1(S)$.

\begin{prop} 
\label{prop:ShearBendGivePleated}
For a given ideal triangulation,
two  pleated surfaces $\bigl(
\widetilde f , r\bigr)$ and
$\bigl( \widetilde f' , r'\bigr)$
with pleating locus $\lambda$ are
isometric if and only if they
have the same exponential
shear-bend factors
$x_i \in \mathbb C^*$ at the
components
$\lambda_i$ of $\lambda$.
Conversely, any set of weights 
$x_i \in \mathbb C^*$ on the
components
$\lambda_i$ of $\lambda$ can be
realized as the exponential
shear-bend parameters of a
pleated surface $\bigl(
\widetilde f , r\bigr)$ with
pleating locus
$\lambda$. \qed
\end{prop}

Note that, for a pleated surface 
$\bigl( \widetilde f , r\bigr)$,
the homomorphism $r\col \pi_1(S)
\rightarrow
\mathrm{PSL}_2(\mathbb C)$ is
completely determined by the map 
$\widetilde f \col
\widetilde S \rightarrow \mathbb
H^3$. The map $\widetilde f$ adds
more data to $r$ as follows. Let
$A \subset S$ be the union of
small annulus neighborhoods of
all the punctures of $S$.  There
is a one-to-one correspondence
between the components of the
preimage $\widetilde A$ of
$A$ and the \emph{peripheral
subgroups} of
$\pi_1(S)$, namely of the images
of the homomorphisms $\pi_1(A)
\rightarrow \pi_1(S)$ defined by
all possible choices of base
points and paths joining these
base points. For a component
$\widetilde A_\pi$ of
$\widetilde A$ corresponding to a
peripheral subgroup $\pi \subset
\pi_1(S)$, the images under
$\widetilde f$ of the triangles of
$\widetilde S - \widetilde
\lambda$ that meet
$\widetilde A_\pi$ all have a
vertex $z_\pi$ in common in
$\widehat {\mathbb C} =
\partial_\infty \mathbb H^3$, and
this vertex is fixed by $r(\pi)$.
Therefore, $\widetilde f$
associates to each peripheral
subgroup
$\pi$ of $\pi_1(S)$ a point
$z_\pi \in
\partial_\infty \mathbb H^3$
which is fixed under $r(\pi)$. In
addition this assignment is
$r$--equivariant in the sense that
$z_{\gamma\pi\gamma^{-1}}=
r(\gamma) z_\pi$ for every
$\gamma \in \pi_1(S)$. 

By definition, an \emph{enhanced
homomorphism} $(r,
\{z_\pi\}_{\pi \in
\Pi})$ of
$\pi_1(S)$ in
$\mathrm{PSL}_2(\mathbb C)$
consists of a group homomorphism
$r\col \pi_1(S) \rightarrow
\mathrm{PSL}_2(\mathbb C)$
together with an
$r$--equivariant assignment of a
fixed point
$z_\pi \in \partial_\infty
\mathbb H^3$ to each peripheral
subgroup $\pi$ of $\pi_1(S)$.
Here $\Pi$ denotes the set of
peripheral subgroups of
$\pi_1(S)$. By abuse of notation,
we will often write $r$ instead
of 
$(r, \{z_\pi\}_{\pi \in
\Pi})$ of
$\pi_1(S)$.

In general, a homomorphism
$r\col \pi_1(S) \rightarrow
\mathrm{PSL}_2(\mathbb C)$ admits
few possible enhancements.
Indeed, if the  peripheral
subgroup
$r(\pi)$ is parabolic, it fixes
only one point in
$\partial_\infty
\mathbb H^3$ and $z_\pi$ is
therefore uniquely determined by
$r$. If $r(\pi)$ is loxodromic or
elliptic, there are exactly two
possible choices for $z_\pi$,
namely  the end points of the
axis of $r(\pi)$; choosing one of
these points as $z_\pi$ therefore
amounts to choosing an orientation
for the axis of
$r(\pi)$. The only case where
there are many possible choices
for $z_\pi$ is when $r(\pi)$ is
the identity, which is highly
non-generic. 

When all the exponential
shear-bend parameters $x_i
\in\mathbb C^*$ are positive
real, there is no bending and the
associated pleated surface
$\widetilde f$ immerses
$\widetilde S$ in a hyperbolic
plane in
$\mathbb H^3$. In particular, the
associated pleated surface 
$\bigl( \widetilde f , r\bigr)$
can be chosen so that the image of
$r$ is contained in the isometry
group
$\mathrm{PSL}_2(\mathbb R)$ of
the hyperbolic plane $\mathbb
H^2$. It can be shown that 
$r\col \pi_1(S) \rightarrow
\mathrm{PSL}_2(\mathbb R)$ is
injective and has discrete image,
and that each  peripheral
subgroup is either parabolic or
loxodromic; see for instance
\cite[\S 3.4]{Th4}.
 In particular, the
enhanced homomorphism $r$ defines
an element of the enhanced
Teichm\"uller space of
$S$, in the terminology of 
\cite{Liu}. The positive real
parameters
$x_i$  are by definition the
exponential shear coordinates for
the enhanced Teichm\"uller space
of $S$.

Given an ideal triangulation
$\lambda$,
Proposition~\ref{prop:ShearBendGivePleated}
and the above observations
associate to a non-quantum
representation $\rho_\lambda^1\col
\mathcal T^1_\lambda
\rightarrow
\mathrm{End} (\mathbb C)$ an
enhanced homomorphism 
$r_\lambda \col
\pi_1(S)
\rightarrow
\mathrm{PSL}_2(\mathbb C)$.  This
correspondence is particularly
well-behaved as we move from one
ideal triangulation to another.

\begin{lem}
\label{lem:NonQuantumGivesHyper}
Let the ideal triangulation
$\lambda'$ be obtained from
$\lambda$ by a re-indexing or a
diagonal exchange, and consider
two non-quantum representations 
$\rho_\lambda^1\col
\mathcal T^1_\lambda
\rightarrow
\mathrm{End} (\mathbb C)$ and
$\rho_{\lambda'}^1\col
\mathcal T^1_{\lambda'}
\rightarrow
\mathrm{End} (\mathbb C)$ such
that $\rho_{\lambda'}^1 =
\rho_\lambda^1 \circ
\Phi_{\lambda\lambda'}^1$. Then
the pleated surfaces
$\bigl(
\widetilde f_\lambda , r_\lambda
\bigr)$ and
$\bigl(
\widetilde f_{\lambda'} ,
r_{\lambda'} \bigr)$ respectively
associated to $\rho_\lambda^1$
and $\rho_{\lambda'}^1$ define
the same enhanced homomorphism
$r_\lambda = r_{\lambda'} \col
\pi_1(S)
\rightarrow
\mathrm{PSL}_2(\mathbb C)$, up to
conjugation by an element of
$\mathrm{PSL}_2(\mathbb C)$. 
\end{lem}

\begin{proof} The property is
immediate when
$\lambda'$ is obtained by
re-indexing the components of
$\lambda$. We can therefore
suppose that $\lambda'$ is
obtained from
$\lambda$ by a diagonal exchange
along the component $\lambda_i$. 

 For a component
$\widetilde\lambda_i$ of the
preimage of $\lambda_i$, consider
as before the left and right
components $T_{\mathrm l}$ and
$T_{\mathrm r}$ of
$\widetilde S -\widetilde\lambda$
that are adjacent to $\lambda_i$,
the end points $z_+$ and $z_-$ of
$\widetilde f_\lambda\bigl(
\widetilde
\lambda_i \bigr)$, and the
remaining vertices $z_{\mathrm
l}$ and $z_{\mathrm r}$ of the
triangles 
$\widetilde f_\lambda(T_{\mathrm
l})$ and 
$\widetilde f_\lambda(T_{\mathrm
r})$. Let $Q \bigl( \widetilde
\lambda_i \bigr) \subset
\widetilde S$ be the open square
$T_{\mathrm l} \cup T_{\mathrm r}
\cup 
\widetilde
\lambda_i$; it admits $\widetilde
\lambda_i$ as a diagonal, but
also a component $\widetilde
\lambda'_i$ of $\widetilde
\lambda'$ as another diagonal. 

Because $\rho_{\lambda'}^1 =
\rho_\lambda^1 \circ
\Phi_{\lambda\lambda'}^1$ is
well-defined, the exponential
shear-bend parameter $x_i \in
\mathbb C^*$ of 
$\bigl(
\widetilde f_\lambda , r_\lambda
\bigr)$ along $\lambda_i$ is
different from $-1$.  This
implies that the points 
$z_{\mathrm l}$ and $z_{\mathrm
r}$ are distinct. We can
therefore modify
$\widetilde f_\lambda$ on $Q
\bigl( \widetilde
\lambda_i \bigr)$ so that it
sends the diagonal $\widetilde
\lambda_i'$ to the geodesic of
$\mathbb H^3$ joining $z_{\mathrm
l}$ to $z_{\mathrm r}$, and the
square $Q
\bigl( \widetilde
\lambda_i \bigr)$ to the union of
the ideal triangles with
respective vertex sets $\{
z_{\mathrm l}, z_{\mathrm r},
z_+\}$ and $\{ z_{\mathrm l},
z_{\mathrm r}, z_-\}$. As
$\widetilde \lambda_i$ ranges
over all the components of the
preimage of $\lambda_i$ in
$\widetilde\lambda$, the
corresponding squares
$Q
\bigl( \widetilde
\lambda_i \bigr)$ are pairwise
disjoint, and we can therefore
perform this operation
equivariantly with respect to
$r_\lambda$. This gives a pleated
surface $\bigl( \widetilde
f_\lambda' , r_\lambda \bigr)$
with pleating locus
$\lambda'$ and with the same
holonomy $r_\lambda  \col
\pi_1(S)
\rightarrow
\mathrm{PSL}_2(\mathbb C)$ as the
original pleated surface
$\bigl(\widetilde f_\lambda ,
r_\lambda \bigr)$. Note that 
$\bigl( \widetilde f_\lambda' ,
r_\lambda
\bigr)$ even has the same
associated enhanced
homomorphism as $\bigl(
\widetilde f_\lambda , r_\lambda
\bigr)$. 

It remains to show that the
exponential shear-bend parameters
of $\bigl( \widetilde f_\lambda'
, r_\lambda \bigr)$ are the
numbers $x_i'\in
\mathbb C^*$ associated to the
non-quantum representation
$\rho_{\lambda'}^1 =
\rho_\lambda^1 \circ
\Phi_{\lambda\lambda'}^1 \col
\mathcal T^1_{\lambda'}
\rightarrow
\mathrm{End} (\mathbb C)$. The
coordinate change isomorphism 
$\Phi_{\lambda\lambda'}^1 \col
\mathcal T^1_{\lambda'}
\rightarrow \mathcal T^1_\lambda
$ was specially designed so that,
when the
$x_i$ are real positive and
correspond to shear coordinates
of the enhanced Teichm\"uller
space, it exactly reflects the
corresponding change of shear
coordinates for the enhanced
Teichm\"uller space; see for
instance \cite{Liu}. The
corresponding combinatorics of
cross-ratios automatically extend
to the complex case, and
guarantees that the non-quantum
representation $
\mathcal T^1_{\lambda'}
\rightarrow
\mathrm{End} (\mathbb C)$ defined
by the $x_i'$ is exactly
$\rho_{\lambda'}^1 =
\rho_\lambda^1 \circ
\Phi_{\lambda\lambda'}^1 $.

As a consequence, $\bigl(
\widetilde f_\lambda' , r_\lambda
\bigr)$ is isometric to  $\bigl(
\widetilde f_{\lambda'} ,
r_{\lambda'}
\bigr)$, which concludes the
proof.
\end{proof}

\begin{prop}
\label{prop:NonQuantumGivesHyper}
Every representation 
$\rho^1 =
\{\rho_\lambda^1\col
\mathcal T^1_\lambda
\rightarrow
\mathrm{End} (\mathbb
C)\}_{\lambda
\in
\Lambda(S)}$ of the non-quantum
polynomial core $\mathcal T_S^1$
uniquely determines an enhanced
homomorphism $r\col \pi_1(S)
\to \mathrm{PSL}_2(\mathbb C)$ such
that, for every ideal triangulation
$\lambda\in
\Lambda(S)$, $r$ is the enhanced
homomorphism associated to the
pleated surface with bending
locus $\lambda$ and with
exponential shear bend parameters
$\rho_\lambda^1(X_i)
\in \mathbb C^*$, for $i=1$,
\dots, $n$.  Conversely, two
representations of $\mathcal
T_S^1$ that induce the same
enhanced homomorphism $r\col \pi_1(S)
\to \mathrm{PSL}_2(\mathbb C)$ must be equal. 
\end{prop}

\begin{proof} The first statement
is an immediate consequence of
Lemma~\ref{lem:NonQuantumGivesHyper}. 

To prove the second statement,
 suppose that the two
representations $\rho$ and
$\rho'$ of $\mathcal T^1_S$
induce the same enhanced
homomorphism,
consisting of a homomorphism
$r\col \pi_1(S)
\rightarrow
\mathrm{PSL}_2(\mathbb C)$ and of
an $r$--equivariant family of
fixed points $z_\pi$ associated
to the peripheral subgroups $\pi$
of $\pi_1(S)$. Let $\bigl
(\widetilde f_\lambda,
r_\lambda\bigr)$ and
$\bigl(
\widetilde f_\lambda',
r_\lambda'\bigr)$ be the two
pleated surfaces with bending
locus
$\lambda$ respectively associated
to $\rho$ and $\rho'$. After
isometries, we can arrange that
$r_\lambda = r_\lambda' =r$. 

Each end of a component
$\widetilde
\lambda_i$ of the preimage
$\widetilde\lambda
\subset \widetilde S$ specifies
two peripheral subgroups $\pi$
and $\pi'$ of $\pi_1(S)$. By
construction $\widetilde
f_\lambda$ and
$\widetilde f_\lambda'$ must both
send
$\widetilde \lambda_i$ to the
geodesic of
$\mathbb H^3$ joining the two
points $z_\pi$ and $z_{\pi'}$.
After a
$\pi_1(S)$--equivariant isotopy of
$\widetilde S$, one can arrange
that 
$\widetilde f_\lambda$ and
$\widetilde f_\lambda'$ coincide
on
$\widetilde\lambda$, and
eventually over all of
$\widetilde S$ by adjustment on
the triangle components of
$\widetilde S-\widetilde
\lambda$. In particular, the two
pleated surfaces $\widetilde
f_\lambda$ and
$\widetilde f_\lambda'$ now
coincide. Since these pleated
surfaces now  have the
same exponential shear-bend
parameters, it follows that
$\rho$ and $\rho'$ coincide on
$\mathcal T_\lambda^1$, and
therefore over all of $\mathcal
T_S^1$.
\end{proof}

By definition, the enhanced
homomorphism
$r \col \pi_1(S) \rightarrow
\mathrm{PSL}_2(\mathbb C)$
provided by Proposition~\ref
{prop:NonQuantumGivesHyper} is
the
\emph{hyperbolic shadow} of the
non-quantum representation
$\rho^1$. In the case
where $\rho^1$ is the non-quantum
shadow of a
representation $\rho$ of the
polynomial core $\mathcal T_S^q$
of the quantum Teichm\"uller
space (for a primitive $N$--th
root $q$ of $(-1)^{N+1}$), we will
also say that
$r$ is the \emph{hyperbolic
shadow} of
$\rho$. 

Not every enhanced homomorphism
from $\pi_1(S)$ to
$\mathrm{PSL}_2(\mathbb C)$ is
associated to a representation of
the polynomial core $\mathcal
T_S^1$ as above. However, many
geometrically interesting ones
are. 

\begin{lem}
\label{lem:FaithfulAlwaysShadow}
Consider an injective homomorphism
$r\col \pi_1(S)
\rightarrow
\mathrm{PSL}_2(\mathbb C)$. Then,
every enhancement of $r$ is
the hyperbolic shadow of  a
representation
$\rho^1$ of the non-quantum
polynomial core $\mathcal T_S^1$.
\end{lem}

\begin{proof} The key property is
that the stabilizer of a point
$z\in \partial_\infty
\mathbb H^3$ in
$\mathrm{PSL}_2(\mathbb C)$ is
solvable, whereas two distinct
peripheral subgroups of
$\pi_1(S)$ generate a free
subgroup of rank 2, which cannot
be contained in a solvable
group. It follows that any
enhancement of $r$ associates
distinct points
$z_\pi$ and $z_{\pi'}$ to
distinct peripheral subgroups
$\pi$ and $\pi'$. 

Let $\lambda$ be an arbitrary
ideal triangulation of $S$, with
preimage
$\widetilde\lambda$ in the
universal covering
$\widetilde S$. The corners of
each component
$T$ of
$\widetilde S-\widetilde
\lambda$ specify three distinct
peripheral subgroups $\pi^T_1$,
$\pi^T_2$ and $\pi^T_3$. We can
then  construct a pleated surface 
$\bigl(
\widetilde f_\lambda , r \bigr)$
with pleating locus $\lambda$,
equivariant with respect to the
given representation $r$, which
sends each component $T$ of
$\widetilde S-\widetilde
\lambda$ to the ideal triangle of
$\mathbb H^3$ with vertices
$z_{\pi^T_1}$,
$z_{\pi^T_2}$, $z_{\pi^T_3}\in
\partial _\infty \mathbb H^3$.
The pleated surface 
$\bigl(
\widetilde f_\lambda , r \bigr)$
defines a representation
$\rho_\lambda^1 \col \mathcal
T^1_\lambda \rightarrow
\mathrm{End} (\mathbb C)$ whose
associated enhanced
homomorphism consists of $r$
and the $z_\pi$. 

As $\lambda$ ranges over all ideal
triangulations, (the proof of)
Lemma~\ref{lem:NonQuantumGivesHyper}
shows that the $\rho_\lambda^1$
fit together to provide a
representation $\rho^1$ of the
polynomial core $\mathcal T^1_S$
whose associated enhanced
representation consists of $r$
and the $z_\pi$.  
\end{proof}

An injective homomorphism  $r\col
\pi_1(S)
\rightarrow
\mathrm{PSL}_2(\mathbb C)$ admits
$2^l$ enhancements, where $l$  is
the number of ends of $S$ whose
image under
$r$ is loxodromic. Combining
Theorem~\ref
{thm:NonQuantumGivesRep},
Proposition~\ref
{prop:NonQuantumGivesHyper} and
Lemma~\ref
{lem:FaithfulAlwaysShadow}
immediately gives:

\begin{thm}
\label{thm:HyperGivesRep} Let $q$
be a primitive $N$--th root of
$(-1)^{N+1}$. Up to isomorphism,
an  injective homomorphism  
$r\col \pi_1(S)
\rightarrow
\mathrm{PSL}_2(\mathbb C)$ is the
hyperbolic shadow of $2^lN^p$ 
if $N$ is odd, and $ 2^{2g+l}N^p$
if $N$ is even,   irreducible
finite-dimensional
representations of the polynomial
core $\mathcal T^q_S$ (where $g$
is the genus of
$S$, $p$ is its number of
punctures, and $l$ is the number
of ends of $S$ whose image under
$r$ is loxodromic). \qed
\end{thm}

\section{Invariants of surface
diffeomorphisms}
\label{sect:Invariants}

Theorem~\ref{thm:HyperGivesRep}
provides a finite-to-one
correspondence between
representations of the polynomial
core $\mathcal T_S^q$ and certain
homomorphisms from $\pi_1(S)$ to
$\mathrm{PSL}_2(\mathbb C)$. We
will take advantage of this
correspondence to construct
interesting representations of
the polynomial core by using
hyperbolic geometry.

Let $\varphi\col S \rightarrow S$
be an orientation-preserving
diffeomorphism of the surface
$S$. If $\lambda$ is an ideal
triangulation of $S$, $\varphi$
induces a natural isomorphism
$\varphi_\lambda^q
\col \mathcal T_\lambda^q 
\rightarrow \mathcal
T_{\varphi(\lambda)}^q$ which, to
the
$i$--th generator $X_i$ of the
Chekhov-Fock algebra
$\mathcal T_\lambda^q$
corresponding to the component
$\lambda_i$ of $\lambda$,
associates the
$i$--th generator $X_i'$ of
$\mathcal T_{\varphi(\lambda)}^q$
corresponding to the component
$\varphi(\lambda_i)$ of
$\varphi(\lambda)$. The existence
of $\varphi$ guarantees that the
$X_i$ and $X_i'$ satisfy the same
relations, so that
$\varphi_\lambda^q$ is a
well-defined algebra isomorphism. 

The isomorphism
$\varphi_\lambda^q$ induces an
isomorphism
$\widehat\varphi_\lambda^q
\col \widehat{\mathcal
T}_\lambda^q 
\rightarrow \widehat{ \mathcal
T}_{\varphi(\lambda)}^q$ between
the corresponding fraction
algebras. As
$\lambda$ ranges over all ideal
triangulations, the
$\widehat\varphi_\lambda^q$ commute
with the coordinate change
isomorphisms
$\Phi_{\lambda\lambda'}^q$, in 
the sense that
$\widehat\varphi_\lambda^q
\circ \Phi_{\lambda\lambda'}^q=
\Phi_{\varphi(\lambda)
\varphi(\lambda' )}^q
\circ
\widehat\varphi_{\lambda'}^q $.
The
$\widehat\varphi_\lambda^q$
consequently define an isomorphism
$\widehat\varphi^q_S$ of the
quantum Teichm\"uller space
$\widehat {\mathcal T}_S^q$. Note
that
$\widehat\varphi^q_S$ sends the
image of
$\mathcal T_\lambda^q$ in
$\widehat {\mathcal T}_S^q$ to
$\mathcal
T_{\varphi(\lambda)}^q$, and
therefore induces an isomorphism
$\varphi_S^q$ of  the polynomial
core
$\mathcal T_S^q$.

In particular, $\varphi$ now acts
on the set $\mathcal R^q$ of
irreducible finite-dimensional
representations of the polynomial
cores $\mathcal T^q_S$ by
associating to  the representation
$\rho =
\{\rho_\lambda\col
\mathcal T^q_\lambda
\rightarrow
\mathrm{End} (V)\}_{\lambda \in
\Lambda(S)}$  the representation
$\rho \circ \varphi^q_S=
\{\rho_{\varphi(\lambda)} \circ
\varphi_\lambda^q\col
\mathcal T^q_\lambda
\rightarrow
\mathrm{End} (V)\}_{\lambda \in
\Lambda(S)}$.

\begin{lem}
\label{lem:HyperShadowAndDiffeo}
If $\rho$ is an irreducible
finite-dimensional representation
of the polynomial core
$\mathcal T_S^q$ and if the
enhanced homomorphism $(r,
\{z_\pi\}_{\pi \in
\Pi})$ is its hyperbolic shadow,
then the hyperbolic shadow of the
representation $\rho
\circ \varphi^q_S$ is equal to
$(r \circ
\varphi*,
\{z_{\varphi^*(\pi)}\}_{\pi \in
\Pi})$, where $\varphi^* \col
\pi_1(S)
\rightarrow \pi_1(S)$ is the
isomorphism induced by the
diffeomorphism $\varphi\col S
\rightarrow S$ for an arbitrary
choice of a path joining the base
point of $S$ to its image under
$\varphi$.
\end{lem}

Note that, up to isometry of
$\mathbb H^3$, the enhanced
representation
$(r \circ
\varphi*,
\{z_{\varphi^*(\pi)}\}_{\pi \in
\Pi})$ is independent of the
choice of path involved in the
definition of $\varphi^*$.

\begin{proof}[Proof of Lemma~\ref
{lem:HyperShadowAndDiffeo}] Let
$\widetilde
\varphi\col \widetilde S
\rightarrow
\widetilde S$ be an arbitrary
lift of
$\varphi$ to the universal cover
$\widetilde S$.   If
$\rho =
\{\rho_\lambda\col
\mathcal T^q_\lambda
\rightarrow
\mathrm{End} (V)\}_{\lambda \in
\Lambda(S)}$ and if  $\bigl(
\widetilde f_\lambda , r_\lambda
\bigr)$ is the pleated surface
with pleating locus $\lambda$
associated to
$\rho_\lambda$, the pleated
surface with pleating locus
$\lambda$ associated to
$\rho_{\varphi(\lambda)} \circ
\varphi^q_\lambda$ is isometric
to $\bigl(
\widetilde f_{\varphi(\lambda)}
\circ
\widetilde\varphi ,
r_{\varphi(\lambda)}
\circ \varphi^*
\bigr)$. The result then
immediately follows from
definitions. 
\end{proof}

We are now ready to use geometric
data to construct special
representations of the polynomial
core. This construction will
require the  diffeomorphism
$\varphi$ to be
\emph{homotopically aperiodic}
(or \emph{homotopically
pseudo-Anosov}) namely such that,
for every
$n>0$ and every non-trivial
$\gamma \in
\pi_1(S)$, 
$\varphi_*^n(\gamma)$ is not
conjugate to $\gamma$ in
$\pi_1(S)$.  The Nielsen-Thurston
classification of surface
diffeomorphisms
\cite{Th3, FLP} asserts that
every isotopy class of surface
diffeomorphism can be uniquely
decomposed into pieces that are
either periodic or homotopically
aperiodic. 

There is another characterization
of homotopically aperiodic
surface diffeomorphisms in terms
of the geometry of their mapping
torus. The \emph{mapping torus}
$M_\varphi$ of the diffeomorphism
$\varphi\col S
\rightarrow S$ is the
3--dimensional manifold quotient
of
$S\times \mathbb R$ by the free
action of $\mathbb Z$ defined by
$n \cdot (x,t) = (\varphi^n(x),
t+n)$ for $n \in \mathbb Z$ and
$(x,t) \in S\times\mathbb R$.
Thurston's Hyperbolization Theorem
\cite{Th2}  asserts
that
$\varphi$ is homotopically
aperiodic if and only if the
mapping torus
$M_\varphi$ admits a complete
hyperbolic metric; see \cite{Ota}
for a proof of this statement.
When this hyperbolic metric
exists, it is unique by Mostow's
Rigidity Theorem \cite{Mos}, and
its holonomy associates to
$\varphi$ an injective
homomorphism
$r_\varphi\col \pi_1(M_\varphi)
\rightarrow
\mathrm{PSL}_2(\mathbb C)$,
uniquely defined up to
conjugation by an element of 
$\mathrm{PSL}_2(\mathbb C)$,  for
which every peripheral subgroup is
parabolic. Consider the map
$f\col S \rightarrow M_\varphi$
composition of the natural
identification $S = S\times
\{0\} \subset S \times \mathbb R$
and of the projection $S
\times \mathbb R \rightarrow
M_\varphi = S \times \mathbb
R/\mathbb Z$. For a suitable
choice of base points, this
enables us to specify a
restriction $r_\varphi\col
\pi_1(S)
\rightarrow
\mathrm{PSL}_2(\mathbb C)$ of the
holonomy homomorphism of
$M_\varphi$. 

The key property is now that $f$ 
is homotopic to $f \circ \varphi$
in $M_\varphi$. This has the
following immediate consequence.

\begin{lem}
\label{lem:InvariantHypRep} The
homomorphisms $r_\varphi$ and
$r_\varphi \circ \varphi^* \col
\pi_1(S)
\rightarrow
\mathrm{PSL}_2(\mathbb C)$ are
conjugate by an element of
$\mathrm{PSL}_2(\mathbb C)$. \qed
\end{lem}

Since every peripheral subgroup of
$\pi_1(S)$ is parabolic for
$r_\varphi$, the homomorphism
$r_\varphi$ admits a unique
enhancement.
Let $\mathcal
R_\varphi^q \subset \mathcal R^q$
be the set of (isomorphism
classes) of irreducible
finite-dimensional
representations of the polynomial
core $\mathcal T_S^q$ whose
hyperbolic shadow is equal to
$r_\varphi$. By Theorem~\ref
{thm:HyperGivesRep},
 the
set
$\mathcal R_\varphi^q$ is finite,
and has $N^p$ or $2^{2g}N^p$
elements according to whether $N$
is odd or even.  By Lemmas~\ref 
{lem:HyperShadowAndDiffeo}  and
\ref{lem:InvariantHypRep}, the
set $\mathcal R_\varphi^q$ is
invariant under the action of
$\varphi$.

By finiteness of $\mathcal
R_\varphi^q$, for every
$\rho =
\{\rho_\lambda\col
\mathcal T^q_\lambda
\rightarrow
\mathrm{End} (V)\}_{\lambda \in
\Lambda(S)}$ in
$\mathcal R_\varphi^q$, there is a
smallest integer
$k\geq1$  such that
$\rho
\circ
(\varphi_S^q)^k=\rho$ in
$\mathcal R_\varphi^q$. This does
not mean that the
representations $\rho
\circ
(\varphi_S^q)^k$ and $\rho$ of
the polynomial core $\mathcal
T_S^q$ over $V$ coincide,
but that there exists an
automorphism
$L_\rho$  of $V$ such that
\begin{equation*}
\rho
\circ
(\varphi_S^q)^k (X) =
L_\rho \cdot \rho (X)
\cdot L_\rho^{-1}
\end{equation*} 
in
$\mathrm{End}(V)$ for every
$X\in \mathcal T_S^q$, if
we denote by $\cdot$ the
composition in $\mathrm {End}(V)$
and by $\circ$ any other
composition of maps to avoid
confusion. 

\begin{prop} 
\label{prop:DiffeoInvariantCycle}
The automorphism
$L_\rho$ of $V$ depends uniquely
on the orbit of $\rho \in \mathcal
R_\varphi^q$ under $\varphi_S^q$,
up to conjugation by an
automorphism of
$V$ and scalar multiplication by a
non-zero complex number.
\end{prop}

\begin{proof} By
irreducibility of
$\rho$, the
isomorphism
$L_\rho$ of $V$ is completely
determined up to scalar
multiplication by the property
that 
$\rho
\circ
(\varphi_S^q)^k (X) =
L_\rho \cdot \rho_\lambda (X)
\cdot L_\rho^{-1}$ for every
$X\in \mathcal T_S^q$. It is also
immediate that we can take
$L_{\rho\circ \varphi_S^q} =
L_\rho$.  Finally, one needs to
remember that the representation
$\rho$ was considered up to
isomorphism of representations. A
representation isomorphism
replaces
$L_\rho$ by a conjugate.
\end{proof}

We consequently have associated
to each orbit of the action of
$\varphi$ on 
$\mathcal R^q_\varphi$  a square
matrix
$L_\rho$ of rank $N^{3g+p-3}$ or
$N^{3g+p-3}/2^g$, according to
wether $N$ is odd or even, which
is well-defined up to conjugation
and scalar multiplication. It is
not too hard to determine these
orbits in terms of the action of
$\varphi$ on the punctures of
$S$ and, when $N$ is even, on
$H_1(S;\mathbb Z_2)$. However,
this process can be cumbersome. 

Fortunately, when $N$ is odd,
there is preferred fixed point
for the action of $\varphi$ on
$\mathcal R_\varphi^q$. This is
based on the following geometric
observation. Recall from
Lemma~\ref
{lem:CoordChange&HP} that the
central elements $P_j$ associated
to the punctures of $S$ and the
square root $H$ of $P_1P_2\dots
P_p$ are well-defined elements of
the polynomial core $\mathcal
T_S^q$. 

\begin{lem}
\label{lem:SendTo1}
Let $\rho_\varphi^1$ be the
non-quantum representation of
$\mathcal T_S^1$ whose
hyperbolic shadow is equal to
$r_\varphi$. Then
$\rho_\varphi^1$ sends the
central elements $H$ and $P_j$ to
the identity. 
\end{lem}
\begin{proof}
Fix an ideal triangulation
$\lambda$, and let $\bigl(
\widetilde f, r_\varphi \bigr)$
be the pleated surface with
pleating locus $\lambda$
associated to $r_\varphi$. 

Consider the $j$--th puncture
$v_j$ of $S$. Because the
corresponding peripheral subgroup
of
$\pi_1(S)$ is parabolic for
$r_\varphi$, the product of the
exponential shear-bend 
coordinates $x_i \in \mathbb C^*$
associated to the components
$\lambda_i$ converging towards
$v_j$ (counted with multiplicity)
is equal to 1; see for instance
\cite{Bon}. By definition of
$P_j$, this means that the
representation $\mathcal
T_\lambda^1 \rightarrow
\mathrm{End}(\mathbb C)$
induced by $\rho_\varphi^1$ sends
$P_j$ to the identity. 

Since $H^2 = P_1 P_2 \dots P_p$,
it follows that $\rho_\varphi^1$
sends $H$ to $\pm 1 = \pm
\mathrm{Id}_{\mathbb C}$. By
construction \cite{Ota}, the
homomorphism $r_\varphi$ is in
the same component as the fuchsian
homomorphisms in the space of
injective homomorphisms $r \col
\pi_1(S) \rightarrow
\mathrm{PSL}_2(\mathbb C)$. For a
fuchsian homomorphism, all the
$x_i$ are real positive, so that
$\rho_r^1(H)=+1=
\mathrm{Id}_{\mathbb C}$ for the
associated representation. By
connectedness, it follows that 
$\rho_\varphi^1(H)=+1=
\mathrm{Id}_{\mathbb C}$.
\end{proof}

When $N$ is odd, we can paraphrase
Theorem~\ref{thm:ClassRepCFodd}
by saying that a representation
$\rho$ of the Chekhov-Fock
algebra $\mathcal T_\lambda^q$ is
classified by its non-quantum
shadow $\rho^1
\col\mathcal T^1_\lambda
\rightarrow
\mathrm{End}(\mathbb C) = \mathbb
C^*$ and by the choice of an
$N$--th root for $\rho^1(H)$ and
for each of the $\rho^1(P_j)$. In
the case when
$\rho^1 = \rho_\varphi^1$,
Lemma~\ref{lem:SendTo1} provides
an obvious choice for these
$N$--th roots, namely 1.
Therefore, $r_\varphi$ specifies
a unique representation
$\rho_\varphi$ of the polynomial
core $\mathcal T_S^q$ over a
vector space $V$ of dimension
$N^{3g+p-3}$, for which
$\rho_\varphi(H) =
\rho_\varphi(P_j) =
\mathrm{Id}_V$. We can paraphrase
this last condition by saying that
$\rho_\varphi$ induces a
representation of the quantum
cusped Teichm\"uller space, as
defined in \cite{Liu}. 

Since the
action of $\varphi$ on the
polynomial core $\mathcal
T_S^q$ respects
$H$ and permutes the $P_j$, it
follows that the representation
$\rho_\varphi$ is fixed under the
action of $\varphi$. As above,
this means that there exists an
isomorphism $L_\varphi$ of $V$
such that 
\begin{equation*}
\rho_\varphi \circ
\varphi_S^q (X) =
L_\varphi \cdot \rho_\varphi (X)
\cdot L_\varphi^{-1}
\end{equation*} 
in
$\mathrm{End}(V)$ for every
$X\in \mathcal T_S^q$. 

\begin{thm}
\label{thm:DiffeoInvariant}
Let $N$ be odd, and let $q$ be
a primitive $N$--th root of
$1$. The isomorphism
$L_\varphi$ of $V$ defined
above depends uniquely on $q$ and
on the homotopically aperiodic
diffeomorphism $\varphi\col S
\rightarrow S$, up to conjugation
 and up to scalar multiplication.
\qed
\end{thm}

In particular, any invariant of
$L_\varphi$ is an invariant of
$\varphi$. For instance, we
can consider the spectrum of
$\varphi$ (consisting of $3g+p-3$
non-zero complex numbers) up to
scalar multiplication. Similarly,
we can normalize the matrix
$L_\varphi$ so that its determinant
is equal to $1$; its trace
$\mathrm{Tr}(L_\varphi)$ then is a
weaker invariant well-defined up
to  a root of unity. Another
interesting invariant is
$\mathrm{Tr}(L_\varphi)
\mathrm{Tr}(L_\varphi^{-1})$,
which is the trace of the linear
automorphism of $\mathrm{End}(V)$
defined by conjugation by
$L_\varphi$. 

See \cite{Liu2} for explicit
computations of $L_\varphi$ for
diffeomorphisms of the
once-punctured torus and of the
4--times punctured sphere.

\end{document}